\documentclass[Final]{siamltex}
\usepackage[margin=1.2in]{geometry}
\usepackage{amsfonts}
\usepackage{mathrsfs}
\usepackage{amsmath}
\usepackage{amssymb}
\usepackage{amsmath}
\usepackage{float}
\usepackage{booktabs}
\usepackage{multirow}
\usepackage{ifpdf}
\usepackage{enumerate}
\usepackage{graphicx}
\usepackage{color}
\usepackage{lineno}

\usepackage{caption}
\usepackage[justification=centering]{caption}
\usepackage{bm}

\makeatletter

\newcommand{\Rmnum}[1]{\expandafter\@slowromancap\romannumeral #1@}
\makeatother
\newcommand{\R}{\mathbb{R}}

\newtheorem{algorithm}[theorem]{Algorithm}

\usepackage{xcolor}	 
\definecolor{mygreen}{RGB}{44,85,17}
\definecolor{myblue}{RGB}{34,31,217}
\definecolor{mybrown}{RGB}{194,164,113}
\definecolor{myred}{RGB}{255,66,56}
\definecolor{mypurple}{RGB}{255 250 205}
\definecolor{myjackieblue}{RGB}{11,23,70}
\definecolor{mygrey}{RGB}{230 230 250}
\definecolor{rot}{rgb}{0.000,0.000,0.000}

\definecolor{myblue}{rgb}{0.000,0.000,1.000}

\definecolor{samplecolor}{rgb}{0.8,0.8,0.8}
\newcommand{\tcs}{\textcolor{samplecolor}}
\definecolor{hatcolor}{rgb}{0.6,0.6,0.6}
\newcommand{\tch}{\textcolor{hatcolor}}
\definecolor{numcolor}{rgb}{0.1,0.1,0.1}
\newcommand{\tcn}{\textcolor{numcolor}}

\begin{document}

\title{Bayesian Approach to Inverse Time-harmonic Acoustic Scattering from Sound-soft Obstacles with Phaseless Data}


\author{Zhipeng Yang\footnotemark[3]\ , Xinping Gui\footnotemark[3]\ , Ju Ming\footnotemark[2] \and Guanghui Hu\footnotemark[1]\ }

\date{}
\maketitle

\renewcommand{\thefootnote}{\fnsymbol{footnote}}
\footnotetext[2]{School of Mathematics and Statistics, Huazhong University of Science and Technology, Wuhan 430074, P. R. China, (jming@hust.edu.cn).}
\footnotetext[3]{Department of Applied Mathematics, Beijing Computational Science Research Center, Beijing 100193, P. R. China, (yangzhp@csrc.ac.cn), (gui@csrc.ac.cn).}
\footnotetext[1]{School of Mathematical Sciences, Nankai University, Tianjin 300071, P. R. China, (ghhu@nankai.edu.cn, corresponding author).}

\renewcommand{\thefootnote}{\arabic{footnote}}

\begin{abstract}
This paper concerns the Bayesian approach to inverse acoustic scattering problems of inferring the position and shape of a sound-soft obstacle from phaseless far-field data generated by point source waves.
To improve the convergence rate, we use
the Gibbs sampling and preconditioned Crank-Nicolson (pCN) algorithm with random proposal variance to implement the Markov chain Monte Carlo (MCMC) method. This usually leads to heavy computational cost, since the unknown obstacle is parameterized in high dimensions. To overcome this challenge, we examine a surrogate model constructed by the generalized polynomial chaos (gPC) method to reduce the computational cost.
Numerical examples are provided to illustrate the effectiveness of the proposed method.
\end{abstract}

\begin{keywords}
Inverse scattering, phaseless far-field data, Bayesian inference, MCMC, Helmholtz equation
\end{keywords}

\begin{AMS}
    35R30, 35P25, 62F15, 78A46
\end{AMS}

\section{Introduction}
In this paper, we consider inverse scattering problems of reconstructing an acoustically sound-soft obstacle from phaseless far-field data in two dimensions.
Let $D\subset \mathbb{R}^2$ be a sound-soft scatterer, which occupies a bounded subset with $C^{2}$-smooth boundary $\partial D$ such that its exterior $\mathbb{R}^{2} \backslash \overline{D}$ is connected.
Suppose that $D$ is embedded in a homogeneous isotropic background medium and that it is
illuminated by a given time-harmonic incident field $u^{\text{in}}$,  which satisfies the Helmholtz equation
\[
\Delta u^{\rm in}+k^2 u^{\rm in}=0
\]
at least in a neighboring area of $D$. Here $k>0$ denotes the wavenumber of the background medium.
The total field $u$ is defined as $u = u^{\text{in}} + u^{\text{sc}}$ in $\mathbb{R}^{2}\backslash \overline{D}$, where
 $u^{\text{sc}}$ is the corresponding scattered (perturbed) field.
Since $D$ is physically sound-soft, the total field $u$ satisfies
  the Dirichlet boundary condition $u = 0$ on the boundary $\partial D$ together with the Helmholtz equation $\Delta u+k^2 u=0$ in $\R^2\backslash\overline{D}$.
Furthermore, the scattered field $u^{\text{sc}}$ is required to fulfill the Sommerfeld radiation condition (see (\ref{SommerfeldCond}) below) at infinity, leading to
the asymptotic behaviour
\begin{equation}
    u^{\text{sc}}(x)
    = \frac{ e^{\mathrm{i}k|x|} }{ \sqrt{|x|} }
      \left\{ u^{\infty}(\hat{\mathbf{x}})
              + \mathcal{O}\left(\frac{ 1 }{ \sqrt{|x|} }\right)
      \right\},
      \quad |x| \rightarrow \infty,
          \label{asymptotic_behaviour}
\end{equation}
uniformly in all directions $\hat{\mathbf{x}}:=x/|x|\in\mathbb{S}:= \{x: |x|=1 \}$. Here, $i = \sqrt{-1}$ is the imaginary unit and $|\cdot|$ denotes the norm of a vector or the modules of a complex number.
The far-field pattern $u^{\infty}: \mathbb{S}\rightarrow \mathbb{C}$ is known as a  real-analytic function with phase information.
The above model also appears in the TE polarization of time-harmonic electromagnetic scattering from infinitely long and perfectly conducting cylinders.
It follows from \cite[Chapter 3.2]{DColton_2013_InverseAcousticScattering} that the forward scattering problem always admits a unique solution $u\in C^2(\R^2\backslash\overline{D})\cap C(\overline{D})$.

Uniqueness, stability and algorithm to inverse time-harmonic obstacle scattering from phased far-field patterns have been extensively studied; we refer to monographs \cite{CC2014, DColton_2013_InverseAcousticScattering, Akirsch_2008, Akirsch_2011_MathematicalTheoryIOP, GNakamura_RPotthast_2015_InverseModeling} for historical remarks, an overview of recent progresses and the comparison between different approaches.
In a variety of practical applications, the accurate phase of the far-field pattern is usually difficult and expensive to obtain, or even cannot be obtained. For instance, in optics it is not trivial to measure the phase of electromagnetic waves incited at high frequencies.
Instead, the modulus or intensity of the far-field pattern is much easier to achieve.
Hence, we are interested in the inverse scattering problem from the phaseless far-field pattern $|u^{\infty}|$.
If the phase information is absent, the key challenge lies in the translation invariance property (see e.g. \cite{RKress_1997_ScatterCrackModulus}) of the phaseless far-field pattern for incident plane waves.
To broke the translation invariance property, one approach was recently proposed in \cite{BZhang_2018_UniquenessScatterPhaseless, BZhang_HZhang_2018_FastImagingPhaselessFixedFrequency} by using infinitely many sets of superposition of two plane waves with different directions as incoming waves, in which both uniqueness and algorithm were investigated.
 The authors in \cite{DZhang_YGuo_2018_UniquePhaselessBall} made use of incident waves generated by superposition of a fixed plane wave and some point sources.
Based on the idea of  \cite{BZhang_2018_UniquenessScatterPhaseless}, in our previous work \cite{Bayesian_phaseless_disk_line_2020_IP} we adapt the Bayesian approach to the recovery of a sound-soft disk, a line crack and a kite-shaped obstacle with less parameters from phaseless far-field patterns generated by plane waves.

The purpose of this paper is to consider obstacles with complex geometric shapes which can be parameterized in high dimensions, when the incoming waves are excited by the following point source waves
\begin{equation}
    u_{\ell}^{\text{in}}(x) := \frac{i}{4}H^{(1)}_{0}(k |x -  x^{\text{in}}_{\ell} |),
    \quad x \in \mathbb{R}^{2}\backslash\{x^{\text{in}}_{\ell}\},
    \quad \ell = 1, 2, \cdots, L.
    \label{incident_wave_spherical}
\end{equation}
Here, $x^{\text{in}}_{\ell}\in \R^2\backslash\overline{D}$ is the position of the $l$-th point source and $H_{0}^{(1)}(\cdot)$ the Hankel function of the first kind of order zero. All source positions $x^{\text{in}}_{\ell}, \ell = 1, 2, \cdots, L$ are assumed to lie on a large circle $|x|=R$ which contains the underlying obstacle inside. Emphasis of this paper will be placed upon how to reduce computational cost of the Bayesian approach for recovering complex obstacles. It should be remarked that,
the translation invariance property for phaseless far-field pattern generated by plane waves does not apply to point source waves.
To the best of the authors' knowledge, it still remains open the unique determination of a general sound-soft obstacle from phaseless far-field patterns corresponding to the aforementioned incoming waves.
Klibanov proved unique determination of a compactly supported potential of the stationary three-dimensional Schr$\ddot{\text{o}}$dinger equation from the phaseless near-field data incited by an interval of frequencies \cite{MVKlibanov_2014_3DScatterPhaseless}.
This was later extended in \cite{MVKlibanov_2017_3DScatterPhaseless} to the reconstruction of a smooth wave speed in the three-dimensional Helmholtz equation. In a deterministic setting, we refer to
\cite{OIvanyshyn_2007_ShapeResconstruction, OIvanyshyn_2010_3DObstaclePhaseless, RKress_1997_ScatterCrackModulus} for inversion algorithms based on a Newton-type iterative scheme.

Recently, the Bayesian approach have attracted extensive attention
 for inverse problems \cite{BThanh_GOmar_2014_BayesianScattering, Stuart_2010_Bayesian_perspective, JKaipio_ESomersalo_2006_StatisticalInverse,Stuart_2015_Data_Assimilation}.
In \cite{BThanh_GOmar_2014_BayesianScattering, Stuart_2010_Bayesian_perspective}, the authors built up a framework of the well-posedness of the posterior distribution, which was later used in \cite{Stuart_2014_SubsurfaceFlow} to determine the permeability of the subsurface from hydraulic head measurements.
The Bayesian approach has also been used with great success to solve inverse scattering problems with phase far-field data in \cite{ABaussard_2001_BayesianScattringMicrowave,  JGSun_2019_StekloffBayesian, YJWang_2015_BayesianScatteringInterior}.
Following the framework in our previous paper \cite{Bayesian_phaseless_disk_line_2020_IP},  a surrogate model constructed by the generalized polynomial chaos (gPC) will be adopted in this work for recovering complex sound-soft in high dimensions.

The Bayesian method provides a new perspective in the form of statistical inferences to view inverse scattering problems. It could be an alternative approach to inverse scattering, when we come up against challenges from deterministic inversion schemes, such as a good initial guess required in the optimization-based iterative schemes and a large number of observation data in non-iterative sampling methods.
On the other hand,
as a disadvantage, it always requires expensive computational cost due to the following reasons.
(i) Since it is quite difficult or even impossible to gain an analytical form of the posterior distributions, we always choose sampling methods such as the Markov chain Monte Carlo (MCMC) method \cite{Brooks_2011_MCMC, Gamerman_2006_MCMC, Geyer_1992_MCMC} to perform numerical approximation.
However, an accurate estimation of the posterior distribution often requires a sufficient number of samples, especially for high-dimensional unknown parameters.
(ii) In the iteration process of the Markov chain, accepting or rejecting a candidate state usually requires
one or more forward solutions to calculate the associated Hastings ratio. Hence, the Bayesian method
 involves a quite large number of repeated solutions of the forward problem. Consequently,
 the computational cost of the Bayesian inference is prohibitively expensive, especially when the forward problem is computationally intensive.
Unfortunately, in most applications, the forward model is always a nonlinear operator associated with partial differential equations.
Therefore, how to reduce the computational cost of the MCMC method is a key point to implement the Bayesian method.

Roughly speaking, the total computational cost of MCMC is the product of the number of iteration steps and the computational cost of one forward solution, which gives rise to a criterion how to save computational efforts.
In this paper, the preconditioned Crank-Nicolson (pCN) algorithm \cite{Stuart_2013_pCN} is adopted to reduce the number of iteration steps.
Since the unknown parameters are high dimensional, we also adopt the Gibbs sampling \cite{GemanS_GemanS_1984_Gibbs, JSLiu_2008_MC, GNakamura_RPotthast_2015_InverseModeling} to accelerate the convergence of MCMC method.
However, even with these advanced methods the number of iteration steps is still high, which turns out to be the order of magnitude of tens of thousands or even hundreds of thousands. Then we have to resort to the idea of reducing computational cost of the forward scattering model.
Recently, substantial attempts have been made to accelerate the Bayesian method in inverse problems associated with a computationally intensive forward model.
Using piecewise linear interpolation, Ma and Zabaras \cite{XMa_NZabaras_2009_adaptive_SG_BIP} adopted the adaptive hierarchical sparse grid collocation (ASGC) method to construct an approximation of the stochastic forward model.
Similarly, Marzouk and Xiu \cite{YMarzouk_DXiu_2009_SC_BIP} proposed sparse grid stochastic collocation methods to improve the efficiency of Bayesian inference. The latter are based on the generalized polynomial chaos  to construct a stochastic surrogate model of the forward model over the support of the prior distribution.
Yan and Guo \cite{LYan_LGuo_2015_SC_CS_BIP} develop the same idea by combining the sparse grid stochastic collocation method with the compressive sensing (CS) method. They employ the $\ell_{1}$-minimization to construct the stochastic surrogate model.
In \cite{GZhang_DLu_MYe_MGunzburger_CWebster_2013_adaptive_SG_BIP}, the adaptive sparse-grid high-order stochastic collocation (aSG-hSC) method is used to construct the surrogate system of a nonlinear groundwater reactive transport model.
Based on the truncated Karhunen-Lo$\grave{\text{e}}$ve (KL) expansions of the prior distribution, a reduced model \cite{YMMarzouk_HNNajm_2009_gPC_BIP} is constructed by the Galerkin projection onto a polynomial chaos basis.
In \cite{TCui_YMMarzouk_KEWillcox_2015_data_driven_BIP}, a data-driven strategy is employed to construct the reduced-order model by projecting the full forward model onto a reduced subspace. Besides, the surrogate model also characterizes the posterior distribution, since the snapshots of the reduced-order model are adaptively calculated from the posterior distribution during the iterations of MCMC method.
Liao and Li \cite{QLiao_JLai_2019_adaptive_ANOVA_BIP} proposed the Analysis of Variance (ANOVA) method to reduce the forward model both in the statistical space and in the physical space. The reduced basis ANOVA model with respect to the posterior distribution is then used in the MCMC iterations by an adaptive scheme.

The key of the aforementioned methods is to derive a reduced-order and computationally efficient surrogate for high-fidelity, large-scale, computationally costly forward models. This surrogate model is then utilized in place of the original forward model to reduce the computational cost in the MCMC iterations.
In this paper, we employ the generalized polynomial chaos (gPC) method
to construct a stochastic surrogate model.
The forward scattering model is projected onto a limited number of basis functions over the support of the prior distribution. Such a projection approach will be used in our MCMC iterations. Our strategy here is to sacrifice the accuracy of the forward model to get an inexpensive surrogate, especially when the prior is significantly different from the posterior.
The price we pay is to develop a special strategy to connect the MCMC method with this surrogate.

This paper is organized as follows. In section 2, we introduce the deterministic forward scattering problem. Section 3 is devoted to the Bayesian framework to inverse scattering problems with phaseless data. In section 4, we construct the surrogate model for the forward scattering problem. The generalized polynomial chaos method will be adopted to reduce the computational cost of MCMC method. Numerical examples will be reported in section 5 and
 conclusions are finally made in section 6.

\section{Deterministic Forward Scattering Problem}
In this paper we want to recover the position and shape of an unknown sound-soft obstacle from phaseless far-field patterns corresponding to a set of incident point source waves. Before dealing with the inverse problem, we need to formulate the abstract nonlinear operator in the deterministic setting which maps the obstacle parameters to the far-field observation data.

Since the boundary $\partial D\subset \R^2$ is a closed $C^{2}$-smooth curve, we can represent or approximate $\partial D$ by a finite set $\mathbf{Z}$ of variables
\begin{equation}
    \mathbf{Z}:= ( z_{1}, z_{2}, \cdots, z_{N} )^{\top} \in \mathbb{R}^{N},
    \quad N\in \mathbb{N}_0. \label{parameter_obstacle}
\end{equation}
For example, we can use $\mathbf{Z}:= ( a_{1}$, $b_{1}$, $a_{2}$, $b_{2}$, $\cdots$, $a_{N}$, $b_{N} )^{\top}$ to approximate a star-shaped closed curve where $\{(a_n, b_n): n=1,\cdots, N\}$ stands for the Fourier coefficients in the truncated Fourier expansion.
Let $u^{\text{in}}_{\ell}(x), \ell = 1, 2, \cdots, L$ be incident waves given by the formula \eqref{incident_wave_spherical}. The forward scattering problem is to find the scattered field $u^{\text{sc}}$ to the Helmholtz equation
\begin{equation}
    \Delta u^{\text{sc}} + k^{2}u^{\text{sc}} = 0 \hspace{.2 cm} \text{in} \hspace{.2 cm} \mathbb{R}^{2}\backslash \overline{D},
    \label{Helmholtz_eq}
\end{equation}
which satisfies the inhomogeneous Dirichlet boundary condition
\begin{equation}
    u^{\text{sc}} =  -u^{\text{in}} \hspace{.2 cm} \text{on} \hspace{.2 cm} \partial D,
    \label{DiriCond_Helmholtz_eq}
\end{equation}
and the Sommerfeld radiation condition
\begin{equation}
    \lim_{r\rightarrow \infty} \sqrt{r} \left( \frac{\partial u^{\text{sc}}}{\partial r} - iku^{\text{sc}} \right) = 0, \hspace{.2 cm} r = |x|,
    \label{SommerfeldCond}
\end{equation}
uniformly in all directions $\hat{\mathbf{x}} \in \mathbb{S}$.
The far-field pattern of the scattering model \eqref{Helmholtz_eq}-\eqref{SommerfeldCond} with and without phase information can be expressed in terms
 the obstacle parameters $\mathbf{Z}$  by
\begin{equation}
    u^{\infty}(\hat{\mathbf{x}}; \mathbf{Z}, x^{\text{in}}_{\ell}, k),
    \quad
    | u^{\infty}(\hat{\mathbf{x}}; \mathbf{Z}, x^{\text{in}}_{\ell}, k) |,
    \quad
    \ell = 1, 2, \cdots, L,
    \quad \hat{\mathbf{x}} \in \mathbb{S}. \label{far_field_pattern}
\end{equation}
We rewrite the forward scattering problem by the operator $F^{\ell}: \mathbb{R}^{N}\rightarrow C(\mathbb{S})$ as
\begin{equation}
   F^{\ell}(\mathbf{Z}):= | u^{\infty}(\hat{\mathbf{x}}; \mathbf{Z}, x^{\text{in}}_{\ell}, k) |,
    \quad  \hat{\mathbf{x}} \in \mathbb{S},   \label{forward_model}
\end{equation}
which can be regarded as an abstract map from the space of obstacle parameters to the space of phaseless far-field pattern in the continuous sense. From the well-posedness of forward scattering, the operators $F^{\ell }, \ell = 1, 2, \cdots, L$, are continuous but highly nonlinear.

Let $G = (g_{1}, g_{2}, \cdots, g_{M})^{\top}: C(\mathbb{S})\rightarrow \mathbb{R}^M$ be a bounded linear observation operator with $g_{m}: C(\mathbb{S}) \rightarrow \mathbb{R}_{+}$ defined as
\begin{equation}
    g_{m}( |u^{\infty}( \hat{\mathbf{x}} )| )
    :=|u^{\infty}( \hat{\mathbf{x}}_{m} ) |,
    \quad m=1, 2, \cdots, M,
    \label{observation_m}
\end{equation}
where $\{\hat{\mathbf{x}}_{m} \in \mathbb{S} \}_{m=1}^M $ denotes the set of discrete observation directions.
Corresponding to the incident wave $u^{\text{in}}_{\ell}(x)$ and the obstacle parameters $\mathbf{Z} \in \mathbb{R}^{N}$, we denote the map $\mathcal{G}^{\ell}: \mathbb{R}^{N} \rightarrow \mathbb{R}^{M}$ from the obstacle parameter space to observation space as
\begin{equation}
    \begin{split}
        \mathcal{G}^{\ell}(\mathbf{Z})
        &:= G\circ F^{\ell}(\mathbf{Z}) \\
        &= \big(| u^{\infty}(\hat{\mathbf{x}}_{1}; \mathbf{Z}, x^{\text{in}}_{\ell}, k) |,
           | u^{\infty}(\hat{\mathbf{x}}_{2}; \mathbf{Z}, x^{\text{in}}_{\ell}, k) |, \cdots,
           | u^{\infty}(\hat{\mathbf{x}}_{M}; \mathbf{Z}, x^{\text{in}}_{\ell}, k) | \big)^{\top}.
    \end{split}
    \label{forward_observation_l}
\end{equation}
Let $\mathbf{Y}^{\ell} = (y_{1}^{\ell}, y_{2}^{\ell}, \cdots, y_{M}^{\ell})^{\top}$ be the measurement data of the phaseless far-field pattern with the observation noise $\eta^{\ell} = (\eta_{1}^{\ell}, \eta_{2}^{\ell}, \cdots, \eta_{M}^{\ell})^{\top} \in \mathbb{R}^{M}$. Then we can express the observation data as
\begin{equation}
    \mathbf{Y}^{\ell}
    = \mathcal{G}^{\ell}( \mathbf{Z} ) + \eta^{\ell},
    \quad \ell = 1, 2, \cdots, L,
    \label{observation_nosie_l}
\end{equation}
or equivalently,
\begin{equation}
    y_{m}^{\ell}
    = g_{m}( F^{\ell}(\mathbf{Z}) ) + \eta_{m}^{\ell}
    = | u^{\infty}( \hat{\mathbf{x}}_{m}; \mathbf{Z}, x^{\text{in}}_{\ell}, k ) | + \eta_{m}^{\ell},
    \quad m=1, 2, \cdots, M.
    \label{observation_noise_l_m}
\end{equation}
Now our inverse problem can be stated as following: determine the obstacle parameters $\mathbf{Z}\in \mathbb{R}^{N}$ from the observation data $\mathbf{Y}^{\ell} \in \mathbb{R}^{M}$ with the noise pollution $\eta^{\ell} \in \mathbb{R}^{M}$, $\ell = 1, 2, \cdots, L$.

\section{Bayesian Framework}
Within the Bayesian framework, all parameters are random variables and the key issue is to estimate
the posterior distribution of the obstacle parameters  based on the Bayes¡¯ formula \cite{Stuart_2010_Bayesian_perspective} and the given assumptions of the prior distribution and the observation pollution.
Since an explicit expression of the posterior distribution is not available, we adopt the Markov chain Monte Carlo method (MCMC) to get an approximation of the posterior distribution.

\subsection{Posterior distribution}
The prior distribution of the obstacle parameters $\mathbf{Z}$ depends on the distribution of $z_{n}, n=1, 2, \cdots, N$. Let $\{z_{n}\}_{n=1}^N$ be independent Gaussian variables as $z_{n} \sim \mathcal{N}( m_{n}, \sigma_{n} ), n = 1, 2, \cdots, N$, where $m_{n}, \sigma_{n}$ are the mean and variance of the distribution of $z_n$. For simplicity, we assume that $\sigma_{1} = \cdots = \sigma_{N} = \sigma_{pr}$, implying that $ \mathbf{Z} \sim \mathcal{N}( \mathbf{m}_{pr}, \sigma_{pr} \mathbf{I} )$, where $\mathbf{m}_{pr} = (m_{1}, m_{2}, \cdots, m_{N})^{\top} $ and $\mathbf{I}\in \mathbb{R}^{N\times N}$ is the identity matrix. Then the prior distribution $P_{pr}( \mathbf{Z} )$ is given by
\begin{eqnarray}
    P_{pr}(\mathbf{Z}) = (2\pi \sigma_{pr} )^{-\frac{N}{2}}
    \exp \Big( -\frac{1}{2\sigma_{pr}} | \mathbf{Z} - \mathbf{m}_{pr} |^{2} \Big).
    \label{prior_density}
\end{eqnarray}

We assume that the observation pollution $\eta^{\ell}$ is independent of $u^{\infty}$, and drawn from the Gaussian distribution $\mathcal{N}( \mathbf{0}, \Sigma_{\eta}^{\ell} )$ with the density $\rho^{\ell}$, where $\Sigma_{\eta}^{\ell}\in \mathbb{R}^{M\times M}$ is a self-adjoint positive matrix, $\ell = 1, 2, \cdots, L$. From the observation data with noise \eqref{observation_nosie_l}, it follows the relationship $\mathbf{Y}^{\ell} | \mathbf{Z} \sim \mathcal{N}( \mathcal{G}^{\ell}( \mathbf{Z} ), \Sigma_{\eta}^{\ell} ), \ell = 1, 2, \cdots, L$. Define the model-data misfit function $\Phi^{\ell}( \mathbf{Z}; \mathbf{Y}^{\ell} ) : \mathbb{R}^{N} \times \mathbb{R}^{M} \rightarrow \mathbb{R} $ as
\begin{equation}
    \Phi^{\ell} \big( \mathbf{Z}; \mathbf{Y}^{\ell} \big)
    = \frac{1}{2} | \mathbf{Y}^{\ell} - \mathcal{G}^{\ell}( \mathbf{Z} ) |^{2}_{\Sigma_{\eta}^{\ell}},
    \quad \ell = 1, 2, \cdots, L.
    \label{model_data_misfit_l}
\end{equation}
Here the norm $| \cdot |_{\Sigma_{\eta}^{\ell}}$ is defined as
\begin{equation}
    | x |^{2}_{\Sigma_{\eta}^{\ell}}
    := x^{\top} \big( \Sigma_{\eta}^{\ell} \big)^{-1} x,
    \quad x \in \mathbb{R}^{M},
    \quad \ell = 1, 2, \cdots, L.
    \label{norm_data_misfit}
\end{equation}
Hence, the likelihood function $P_{lhd} \big(\mathbf{Y}^{1}, \mathbf{Y}^{2}, \cdots, \mathbf{Y}^{L}; \mathbf{Z} \big)$ is given by
\begin{equation}
    \begin{split}
        P_{lhd} \big( \mathbf{Y}^{1}, \mathbf{Y}^{2}, \cdots, \mathbf{Y}^{L}; \mathbf{Z} \big)
        & =  \prod_{\ell=1}^{L} \rho^{\ell} \big( \mathbf{Y}^{\ell} - \mathcal{G}^{\ell}( \mathbf{Z} ) \big) \\
        & =  (2\pi)^{-\frac{LM}{2}} \prod_{\ell=1}^{L}  \big(  \mbox{det}( \Sigma_{\eta}^{\ell} ) \big)^{ -\frac{1}{2} }
            \exp\Big( -\sum_{\ell=1}^{L} \Phi^{\ell} \big( \mathbf{Z}; \mathbf{Y}^{\ell} \big) \Big).
    \end{split}
    \label{likelihood_L}
\end{equation}
Furthermore, by the Bayes' theorem \cite{Stuart_2015_Data_Assimilation, Stuart_2010_Bayesian_perspective}, the posterior distribution $P_{post}(\mathbf{Z}; \mathbf{Y}^{1}, \mathbf{Y}^{2}, \cdots, \mathbf{Y}^{L})$ is given as
\begin{equation}
   P_{post}(\mathbf{Z}; \mathbf{Y}^{1}, \mathbf{Y}^{2}, \cdots, \mathbf{Y}^{L})
   = C_{z}^{-1}
     \exp\Big( -\sum_{\ell=1}^{L} \Phi^{\ell} \big( \mathbf{Z}; \mathbf{Y}^{\ell} \big)
               -\frac{1}{2\sigma_{pr}} | \mathbf{Z} - \mathbf{m}_{pr} |^{2} \Big).
   \label{posterior_distribution}
\end{equation}
Here $C_{z}$ is the normalization constant
\begin{equation}
   C_{z}
   = \int_{\mathbb{R}^{N}}
     \exp\Big( -\sum_{\ell=1}^{L} \Phi^{\ell} \big( \mathbf{Z}; \mathbf{Y}^{\ell} \big)
               -\frac{1}{2\sigma_{pr}} | \mathbf{Z} - \mathbf{m}_{pr} |^{2} \Big)
     d \mathbf{Z}.
   \label{posterior_normalization}
\end{equation}

The well-posedness arguments of \cite{BThanh_GOmar_2014_BayesianScattering, Stuart_2010_Bayesian_perspective} can be applied to deal with our inverse scattering
problem with the Bayesian approach. We state
the well-posedness of the posterior distribution in the
 theorem below. Its proof relies heavily on the well-posedness of the forward scattering problem, for example, via the integral equation and variational methods. We omit its proof here and refer to
 \cite{ABaussard_2001_BayesianScattringMicrowave,  JGSun_2019_StekloffBayesian, YJWang_2015_BayesianScatteringInterior, Bayesian_phaseless_disk_line_2020_IP} for detailed discussions.
\begin{theorem}\label{Well_Posterior}
    Let $\mu_{pr}$ and $\mu_{post}$ be the probability measures of the prior distribution $P_{pr}$ and the posterior distribution $P_{post}$. Then $\mu_{post}$ is a well-defined probability measure on $\mathbb{R}^{N}$ and absolutely continuous with respect to prior measure $\mu_{pr}$. What's more, the posterior measure $\mu_{post}$ is Lipschitz in the observation data $\big\{ \mathbf{Y}^{\ell} \big\}_{\ell=1}^{L}$, with respect to the Hellinger distance: if $\mu_{post}^{1}$ and $ \mu_{post}^{2} $ are two posterior measures corresponding to data $\big\{ \mathbf{Y}_{1}^{\ell} \big\}_{\ell=1}^{L}$ and $\big\{ \mathbf{Y}_{2}^{\ell} \big\}_{\ell=1}^{L}$, then there exists $C = C(r) > 0$ such that,
    \begin{equation*}
         d_{\text{Hell}}(\mu_{post}^{1}, \mu_{post}^{2})
         \leq C \sum^{L}_{\ell=1} |\mathbf{Y}_{1}^{\ell} - \mathbf{Y}_{2}^{\ell} |,
    \end{equation*}
    for all $\big\{ \mathbf{Y}_{1}^{\ell} \big\}_{\ell=1}^{L}$, $\big\{ \mathbf{Y}_{2}^{\ell} \big\}_{\ell=1}^{L}$ with $ \max\limits_{\ell=1,2,\cdots,L} \big\{ | \mathbf{Y}_{1}^{\ell} |, | \mathbf{Y}_{2}^{\ell} |  \big\} < r $. Here the Hellinger distance is defined by
    \begin{equation}
        d_{\text {Hell}}\left(\mu_{1}, \mu_{2}\right) :=\sqrt{\frac{1}{2} \int\left(\sqrt{\frac{d \mu_{1}}{d \mu_{0}}}-\sqrt{\frac{d \mu_{2}}{d \mu_{0}}}\right)^{2} d \mu_{0}},
        \label{Hellinger_distance}
    \end{equation}
    where $\mu_{1}, \mu_{2}$ are two measures that are absolutely continuous with respect to $\mu_{0}$.
\end{theorem}


\subsection{Markov chain Monte Carlo method}
In the posterior distribution \eqref{posterior_distribution}, it is challenging to give an explicit expression of the normalization constant $C_{z}$ by the integration \eqref{posterior_normalization}. Hence, a suitable numerical method is needed to calculate the posterior distribution. For this purpose we adopt the Markov chain Monte Carlo method (MCMC) \cite{Brooks_2011_MCMC, Gamerman_2006_MCMC, Geyer_1992_MCMC} to generate a large number of samples subject to the posterior distribution. The numerical approximation of the posterior distribution of unknown obstacle parameters can be obtained by statistical analysis on these samples.

In this section we use the Metropolis-Hastings \cite{Hastings_1970, Metropolis_1953} algorithm to construct MCMC samples. Since the dimension $N$ of the space of obstacle parameters is large, the Metropolis-Hastings algorithm may stay at one state for a quite long time with a huge number of iterations.
In each iteration of the Metropolis-Hastings algorithm, we have to choose a proper candidate multi-dimensional sample, which however relies heavily on the dimension of obstacle parameters.
To overcame this challenge we adopt the Gibbs sampling \cite{GemanS_GemanS_1984_Gibbs, JSLiu_2008_MC, GNakamura_RPotthast_2015_InverseModeling}, especially when the dimension is conditioned only on a small number of other dimensions.
In the iteration of the Gibbs sampling, samples can be chosen to be dependent on partial dimensions of obstacle parameters and are not necessarily uniform in all dimensions.
Noticing that $z_{n}$, $n = 1, 2, \cdots, N$ are independent variables, the Gibbs sampling will significantly improve the convergence rate.

As done in our previous work \cite{Bayesian_phaseless_disk_line_2020_IP}, the preconditioned Crank-Nicolson (pCN) algorithm with a random proposal variance \cite{Stuart_2013_pCN} can be applied to generate an ergodic Markov chains and to improve the convergence rate of the MCMC method. Using this scheme, the candidate state of the obstacle parameters $\tilde{\mathbf{Z}}$ can be iteratively updated from the current state (initial guess) $\mathbf{Z}$ through the formula
\begin{equation}
    \tilde{\mathbf{Z}} = \mathbf{m}_{pr} + ( 1-\beta^2 )^{1/2} (\mathbf{Z} - \mathbf{m}_{pr}) + \beta \omega, \label{pCN}
\end{equation}
where $\beta \in [0, 1]$ is the random proposal variance coefficient and $\omega \sim  \mathcal{N}( \mathbf{0},  \Sigma_{pcn}  ) $ is a zero-mean normal random vector with the covariance matrix $\Sigma_{pcn} = \sigma_{pr} \mathbf{I} \in \mathbb{R}^{N\times N}$. The numerical algorithm is described as follows.

\begin{algorithm}\label{Gibbs_Random_Proposal} \textbf{(Gibbs Sampling with Random Proposal Variance)}
\smallskip
\begin{itemize}
  \item Initialize $\mathbf{Z}_{0} \in \mathbb{R}^{N}$ from the prior distribution $P_{pr}(\mathbf{Z})$ and $\beta_{1_{0}} = \beta_{2_{0}} = \cdots = \beta_{N_{0}} \in [0, 1]$.

  \smallskip

  \item Repeat iteration from $\mathbf{Z}_{j} = (z_{1_{j}}, z_{2_{j}}, \cdots, z_{N_{j}})^{\top}$ to $\mathbf{Z}_{j+1} = (z_{1_{j+1}}, z_{2_{j+1}}, \cdots, z_{N_{j+1}})^{\top}$, $j = 0, 1, \cdots, J_{0}-1$. For $n=1, 2, \cdots, N$:

        \begin{enumerate}
            \smallskip
          \item Draw the candidate sate $\tilde{\mathbf{Z}}_{n} = (z_{1_{j+1}}, \cdots, z_{(n-1)_{j+1}}, \tilde{z}_{n_{j}}, z_{(n+1)_{j}}, \cdots, z_{N_{j}})^{\top}$ by modifying the $n$-th component of the current state $\mathbf{Z}_{j, n} = (z_{1_{j+1}}$, $\cdots$, $z_{(n-1)_{j+1}}$, $z_{n_{j}}$, $z_{(n+1)_{j}}$, $\cdots$, $z_{N_{j}})^{\top}$ using the pCN algorithm \eqref{pCN} with the proposal variance coefficient $\beta_{n_{j}}$ as:
               \begin{equation}
                   \tilde{z}_{n_{j}}
                   = m_{n} + ( 1-\beta_{n_{j}}^2 )^{1/2} (z_{n_{j}} - m_{n}) + \beta_{n_{j}} \omega_{n},
                   \quad \omega_{n} \sim \mathcal{N}( 0, \sigma_{pr} ); \label{pCN_RPV}
               \end{equation}

          \item Compute Hasting ratio $\alpha ( \cdot, \cdot ):\mathbb{R}^{N}\times \mathbb{R}^{N}\rightarrow[1,\infty) $ as:
                \begin{equation}
                    \alpha \big( \mathbf{Z}_{j, n}, \tilde{\mathbf{Z}}_{n} \big)
                    = \min\Bigg \{ 1, \exp\bigg(
                                 \sum_{\ell=1}^{L} \Phi^{\ell} \big( \mathbf{Z}_{j, n}; \mathbf{Y}^{\ell} \big)
                                - \sum_{\ell=1}^{L} \Phi^{\ell} \big( \tilde{\mathbf{Z}}_{n}; \mathbf{Y}^{\ell} \big)
                               \bigg)
                          \Bigg \};
                    \label{alpha_Z}
                \end{equation}

          \item Accept or Reject $\tilde{\mathbf{Z}}$: draw $U \sim \mathcal{U}(0,1)$ and then update $\mathbf{Z}_{j, n}$ by the criterion
                \begin{equation}
                    \mathbf{Z}_{j, n+1} =
                                \left\{
                                        \begin{array}{ll}
                                                \tilde{\mathbf{Z}}_{n}, & \hbox{if\; }U \leq \alpha \big( \tilde{\mathbf{Z}}_{n}, \mathbf{Z}_{j, n} \big)  \\
                                                \mathbf{Z}_{j, n}, & \hbox{if\, otherwise;}
                                        \end{array}
                                \right. \label{update_Z}
                \end{equation}

          \item Generate new proposal variance coefficient $\beta_{n_{j+1}}$ from $\beta_{n_{j}}$. First we set
                \begin{equation}
                   \beta_{new} = ( 1-\gamma^2 )^{1/2} \beta_{n_{j}} + \gamma ( \omega_{\beta} - 0.5 ), \hspace{.3cm} \omega_{\beta} \sim \mathcal{U}(0,1), \label{pCN_beta}
                \end{equation}
                with $\gamma \in [0, 1]$. In our case we choose $\gamma = 0.1$. Then $\beta_{n_{j+1}}$ can be updated by
                \begin{equation}
                    \beta_{n_{j+1}} =
                                \left\{
                                        \begin{array}{ll}
                                                \beta_{new}, & \mbox{if}\quad \beta_{new} \in [0, 1], \\
                                                - \beta_{new}, & \mbox{if}\quad \beta_{new} < 0, \\
                                                \beta_{new} - 1, & \mbox{if}\quad\beta_{new} > 1 .
                                        \end{array}
                                \right. \label{update_beta}
                \end{equation}

            \item Set $\tilde{\mathbf{Z}}_{1} = ( \tilde{z}_{1_{j}}, z_{2_{j}}, \cdots, z_{N_{j}})^{\top}$, $\tilde{\mathbf{Z}}_{N} = (z_{1_{j+1}}, z_{2_{j}}, \cdots, z_{(N-1)_{j+1}}, \tilde{z}_{N_{j}})^{\top}$, $\mathbf{Z}_{j, 1} = \mathbf{Z}_{j}$, $\mathbf{Z}_{j, N} = (z_{1_{j+1}}, z_{2_{j+1}}, \cdots, z_{(N-1)_{j+1}}, z_{N_{j}})^{\top}$, $ \mathbf{Z}_{j+1} = \mathbf{Z}_{j, N+1}$.
        \end{enumerate}

        \smallskip

    \item Select $\mathbf{Z}_{\tilde{j}}$ with the indices $ \tilde{j} = J_{1} + (\hat{j} - 1 ) J_{2}, \   \hat{j} = 1, 2, \cdots, J_{3}$.
\end{itemize}
\end{algorithm}

In the Algorithm \ref{Gibbs_Random_Proposal}, the number $J_{0}, J_{1}, J_{2}, J_{3} \in \mathbb{N}_0$ are four positive integers. The integer $J_{0}$ is the number of  total iterations and $J_{1}$  the number of initial states which will be threw away to ensure that the Markov chain converges to the posterior distribution. By the number $J_{2}$ we mean that $J_{2}$ sates are taken to guarantee the independence of the selected sates. The integer $J_{3}$ denotes the number of totally selected states to approximate the posterior distribution.

\section{Stochastic Surrogate Model}
The computational cost of each iteration in the Algorithm \ref{Gibbs_Random_Proposal} is dominated by the computational cost of the forward map $\mathcal{G}^{\ell}(\mathbf{Z})$, $\ell = 1, 2, \cdots, L$. The total computational cost $T$ of the Algorithm \ref{Gibbs_Random_Proposal} is given by
\begin{equation}
    T = T_{0}NJ_{0},
    \label{cost_Gibbs_01}
\end{equation}
where $T_{0} $ denotes the computational cost of the forward maps $\mathcal{G}^{\ell}(\mathbf{Z}), \ell = 1, 2, \cdots, L$, for $L$ incident waves.
In the previous section the Gibbs sampling and pCN algorithm with a random proposal variance have been adopted to reduce the number of the total iterations $J_{0}$. Below we discuss how to reduce the computational cost $T_{0}$.
Our idea is to adopt the generalized polynomial chaos method (gPC) to construct a surrogate model for the forward map $\mathcal{G}^{\ell}(\mathbf{Z}), \ell = 1, 2, \cdots, L$.

Let $x \sim \mathcal{N}(0, 1)$ be the one-dimensional standard Gaussian variable with the distribution $P_{0}(x) = \frac{1}{\sqrt{2\pi}} e^{-x^2/2}$. The corresponding Gaussian space $\mathcal{L}^{2}_{P_{0}}(\mathbb{R})$ is defined as
\begin{equation}
    \mathcal{L}^{2}_{P_{0}}(\mathbb{R}) := \Big\{
    f(x): \int_{\mathbb{R}} P_{0}(x) | f(x) |^{2} dx < \infty
    \Big\}.
    \label{L2_Gauss_space_stand}
\end{equation}
The normalized one-dimensional Hermite polynomials of order $m$ are defined as
\begin{equation}
    h_{m}(x) := (-1)^{m} e^{x^{2}/2} \frac{d^{m}}{dx^{m}} e^{-x^{2}/2},
    \quad m = 0, 1, 2, \cdots,
    \quad x \in \mathbb{R}.
    \label{1D_Hermite}
\end{equation}
It is well known that the set $\{ h_{m}(x) \}_{m=0}^{\infty}$ is a complete orthonormal basis of $\mathcal{L}^{2}_{P_{0}}(\mathbb{R})$ with respect to the Gaussian distribution $P_{0}(x)$, that is,
\begin{equation}
    \mathbb{E}[ h_{n} h_{m} ]
    = \int_{\mathbb{R}} P_{0}(x)  h_{n}(x)  h_{m}(x) dx
    = \delta_{nm},
    \quad n, m = 0, 1, 2, \cdots.
    \label{orthonormal_1D_Hermite}
\end{equation}
Here $\delta_{nm}$ is the Kronecker delta function.

Recall that the components of the obstacle parameters $\mathbf{Z}$ are Gaussian variables with the prior distribution $P_{pr, n}(z_{n})$ given by
\begin{equation}
    z_{n} \sim \mathcal{N}( m_{n}, \sigma_{pr} ),
    \quad P_{pr, n}(z_{n}) = \frac{1}{\sqrt{2\pi \sigma_{pr}}}
          \exp\Big(-\frac{1}{2\sigma_{pr}} (z_{n} - m_{n})^2\Big),
    \quad n = 1, 2, \cdots, N.
    \label{prior_zn}
\end{equation}
Then for each $z_{n}$, the Gaussian space $\mathcal{L}^{2}_{P_{pr, n}}(\mathbb{R})$ has a complete orthonormal basis $\{ h_{n, m}(x) \}_{m=0}^{\infty}$, which can be obtained by modifying the set $\{ h_{m}(x) \}_{m=0}^{\infty}$ given by (\ref{1D_Hermite}).
 Clearly, the tensor product of the elements of $\{ h_{n, m}(x) \}_{m=0}^{\infty}$, $n = 1, 2, \cdots, N$ form a complete basis of the corresponding $N$-dimensional Gaussian probability space $\mathcal{L}^{2}_{P_{pr}}(\mathbb{R}^{N})$ with respect to the prior distribution $P_{pr}( \mathbf{Z} )$.

Let $\mathcal{I}$ denote the $N$ dimensional multi-indexes:
\begin{equation}
    \mathcal{I} := \{ \alpha = (\alpha_{1}, \alpha_{2}, \cdots, \alpha_{N} )
        : \  \alpha_{n} \in \mathbb{N}_{0},
          \ | \alpha | = \sum_{n=1}^{N} \alpha_{n} < \infty  \}.
    \label{ND_index}
\end{equation}
Then a complete basis of the space $\mathcal{L}^{2}_{P_{pr}}(\mathbb{R}^{N})$ is given by the set $\{ \mathcal{H}_{\alpha}(\mathbf{Z}) \}_{\alpha \in \mathcal{I}}$, defined by
\begin{equation}
    \mathcal{H}_{\alpha}(\mathbf{Z})
    := \prod_{n=1}^{N} h_{\alpha_{n}}(z_{n}),
    \quad \alpha \in \mathcal{I}.
    \label{basis_ND}
\end{equation}
By properties of the operator $\mathcal{G}^{\ell}$, we claim that $| u^{\infty}(\hat{\mathbf{x}}_{m}; \mathbf{Z}, x^{\text{in}}_{\ell}, k) | \in \mathcal{L}^{2}_{P_{pr}}(\mathbb{R}^{N})$, $\ell = 1, 2, ..., L$, $m = 1, 2, ..., M$. By the Cameron-Martin theorem \cite{RHCameron_WTMartin_1947_gPC}, the phaseless data
 $| u^{\infty}(\hat{\mathbf{x}}_{m}; \mathbf{Z}, x^{\text{in}}_{\ell}, k) |$ can be expanded into the series
\begin{equation}
    | u^{\infty}(\hat{\mathbf{x}}_{m}; \mathbf{Z}, x^{\text{in}}_{\ell}, k) |
    := \sum_{\alpha \in \mathcal{I}}
       u_{\alpha}^{\ell, m}\;  \mathcal{H}_{\alpha}(\mathbf{Z}),
    \quad \ell = 1, 2, ..., L,
    \quad m = 1, 2, ..., M,
    \label{u_gPC_infty}
\end{equation}
where $u_{\alpha}^{\ell, m}\in \mathbb{C}$ are referred to as the chaos coefficients given by
\begin{equation}
    u_{\alpha}^{\ell, m}
    = \mathbb{E}\big[
        | u^{\infty}(\hat{\mathbf{x}}_{m}; \mathbf{Z}, x^{\text{in}}_{\ell}, k) |
        \mathcal{H}_{\alpha}(\mathbf{Z}) \big]
    = \int_{\mathbb{R}^{N}}
        P_{pr}(\mathbf{Z})
        | u^{\infty}(\hat{\mathbf{x}}_{m}; \mathbf{Z}, x^{\text{in}}_{\ell}, k) |
        \mathcal{H}_{\alpha}(\mathbf{Z})
        d \mathbf{Z}.
    \label{u_alpha_l_m}
\end{equation}
In this paper, we define the surrogate $\tilde{\mathcal{G}}^{\ell}$ of the forward operator $\mathcal{G}^{\ell}$, $\ell = 1, 2, \cdots, L$ through the gPC approximation of order $\tilde{N}\in \mathbb{N}_{0}$ as
\begin{equation}
    \tilde{\mathcal{G}}^{\ell}(\mathbf{Z})
    := \Big(
          | \tilde{u}^{\infty}(\hat{\mathbf{x}}_{1}; \mathbf{Z}, x^{\text{in}}_{\ell}, k) |,
        \ | \tilde{u}^{\infty}(\hat{\mathbf{x}}_{2}; \mathbf{Z}, x^{\text{in}}_{\ell}, k) |,
        \ \cdots,
        \ | \tilde{u}^{\infty}(\hat{\mathbf{x}}_{M}; \mathbf{Z}, x^{\text{in}}_{\ell}, k) |
    \Big)^{\top},
    \label{forward_observation_surrogate}
\end{equation}
where
\begin{equation}
    | \tilde{u}^{\infty}(\hat{\mathbf{x}}_{m}; \mathbf{Z}, x^{\text{in}}_{\ell}, k) |
    = \sum_{ \alpha \in \mathcal{I}, |\alpha| = 0}^{|\alpha| = \tilde{N}}
       u_{\alpha}^{\ell, m}  \mathcal{H}_{\alpha}(\mathbf{Z}),
    \quad \ell = 1, 2, ..., L,
    \quad m = 1, 2, ..., M.
    \label{u_gPC_finite}
\end{equation}
We adopt the Monte Carlo method \cite{Robert_2004_MC} to calculate the chaos coefficients $u_{\alpha}^{\ell, m}$ in \eqref{u_alpha_l_m}. More details will be presented in  Example 4 of the subsequent section.
The surrogate model $\big\{ \tilde{\mathcal{G}}^{\ell} \big\}_{\ell=1}^{L}$ given by \eqref{forward_observation_surrogate} will be used in place of the original operator $\big\{ \mathcal{G}^{\ell} \big\}_{\ell=1}^{L}$ to generate the candidate sate $\tilde{\mathbf{Z}}_{n}$, $n = 1, 2, \cdots, N$ (see \eqref{pCN_RPV}) in the Algorithm \ref{Gibbs_Random_Proposal}. We summarize the algorithm as follows.

\begin{algorithm}\label{Gibbs_gPC}\textbf{(Gibbs Sampling with Surrogate Model)}
Let $j = 0, 1, \cdots, J_{0} - 1$, $n = 1, 2, \cdots, N$ be given in the first step of Algorithm \ref{Gibbs_Random_Proposal}. We replace the formula \eqref{pCN_RPV} for  generating a candidate state $\tilde{z}_{n_{j}}$ by the following steps:
    \begin{enumerate}
        \item Drew $\hat{J}_{1}$ candidate sates $\tilde{\mathbf{Z}_{n}^{\hat{j}}}$, $\hat{j}=1, 2, \cdots, \hat{J}_{1}$, by the pCN algorithm \eqref{pCN} with the proposal variance coefficient $\beta_{n_{j}}$.
        Here the $n$-th component of $\tilde{\mathbf{Z}}_{n}^{\hat{j}}$ is given by:
               \begin{equation}
                   \tilde{z}_{n_{j}}^{\hat{j}}
                   = ( 1-\beta_{n_{j}}^2 )^{1/2} z_{n_{j}} + \beta_{n_{j}} \omega_{n}^{\hat{j}},
                   \quad \omega_{n}^{\hat{j}} \sim \mathcal{N}( 0, \sigma_{pr} ),
                   \quad \hat{j}=1, 2, \cdots, \hat{J}_{1};
                   \label{pCN_RPV_multi}
               \end{equation}

        \item Compute $\tilde{\phi}^{\hat{j}} = \sum\limits_{\ell=1}^{L} \tilde{\Phi}^{\ell} \big( \tilde{\mathbf{Z}}_{n}^{\hat{j}}; \mathbf{Y}^{\ell} \big) $, $\hat{j}=1, 2, \cdots, \hat{J}_{1}$. Here the function $\tilde{\Phi}^{\ell}$ is defined by (cf. the model-data function \eqref{model_data_misfit_l})
            \begin{equation}
                \tilde{\Phi}^{\ell} \big( \mathbf{Z}; \mathbf{Y}^{\ell} \big)
                = \frac{1}{2} | \mathbf{Y}^{\ell} - \tilde{\mathcal{G}}^{\ell}( \mathbf{Z} ) |^{2}_{\Sigma_{\eta}^{\ell}},
                \quad \ell = 1, 2, \cdots, L.
                \label{model_data_misfit_l_gPC}
            \end{equation}
            Noting that the forward operator $\mathcal{G}^{\ell}$ in \eqref{model_data_misfit_l} has been substituted by the surrogate $\tilde{\mathcal{G}}^{\ell}$ in (\ref{model_data_misfit_l_gPC});

        \smallskip

        \item Sort the candidate sate $ \big\{ \tilde{\mathbf{Z}}_{n}^{\hat{j}} \big \}_{\hat{j}=1 }^{\hat{J}_{1}}$ that corresponds to the sequence $\big\{ \tilde{\phi}^{\hat{j}} \big\}_{\hat{j} =1}^{\hat{J}_{1}}$ from smallest to largest;

        \smallskip

        \item Choose the first $\hat{J}_{2}$ sates of the sorted candidate sates $ \big\{ \tilde{\mathbf{Z}}_{n}^{\hat{j}} \big \}_{\hat{j}=1 }^{\hat{J}_{1}}$. Then compute $\phi^{\hat{j}} = $ \\ $ \sum\limits_{\ell=1}^{L} \Phi^{\ell} \big( \tilde{\mathbf{Z}}_{n}^{\hat{j}}; \mathbf{Y}^{\ell} \big) $, $\hat{j}=1, 2, \cdots, \hat{J}_{2}$, with the model-data function $\Phi^{\ell}$;

        \smallskip

        \item Sort the candidate sate $ \big\{ \tilde{\mathbf{Z}}_{n}^{\hat{j}} \big \}_{\hat{j}=1 }^{\hat{J}_{2}}$
            that corresponds to the sequence $\big\{ \phi^{\hat{j}} \big\}_{\hat{j} =1}^{\hat{J}_{2}}$ from smallest to largest;

        \smallskip

        \item Set the candidate sate $\tilde{\mathbf{Z}}_{n} $ as the first sate of the sorted candidate sates $ \big\{ \tilde{\mathbf{Z}}_{n}^{\hat{j}} \big \}_{\hat{j}=1 }^{\hat{J}_{2}}$.

    \end{enumerate}
\end{algorithm}
We remark that, with the surrogate model the new Algorithm \ref{Gibbs_gPC} explores $\hat{J}_{1}$ states for each iteration in the MCMC method, whereas the Algorithm \ref{Gibbs_Random_Proposal} explores one state only.
Let $T_{1}$ be the computation cost in calculating $\big\{ \mathcal{H}_{\alpha}(\mathbf{Z}^{\hat{j}}) \big\}_{\alpha \in \mathcal{I}, |\alpha| = 0}^{ |\alpha| = \tilde{N}}$, $\hat{j}=1, 2, \cdots, \hat{J}_{1}$.
The computational cost of the Algorithm \ref{Gibbs_gPC} is
\begin{equation}
    \hat{T} = (T_{1} + T_{0}\hat{J}_{2})NJ_{0}.
    \label{cost_Gibbs_gPC}
\end{equation}

Assume that there are also $\hat{J}_{1}$ states at each iteration
in the Algorithm \ref{Gibbs_Random_Proposal}. Without using the surrogate model $\big\{ \tilde{\mathcal{G}}^{\ell} \big\}_{\ell=1}^{L}$, these $\hat{J}_{1}$ states can be evaluated by the model-data function $\Phi^{\ell}$ given by \eqref{model_data_misfit_l}.
To further compare the previous two schemes, we rewrite the Algorithm \ref{Gibbs_Random_Proposal} in the form of Algorithm \ref{Gibbs_gPC} as follows.

\begin{algorithm}\label{Gibbs_Random_Proposal_multi} \textbf{(Gibbs Sampling with Multi Candidate)}
    \begin{enumerate}
        \item Drew $\hat{J}_{1}$ candidate sates $\tilde{\mathbf{Z}}_{n}^{\hat{j}}$, $\hat{j}=1, 2, \cdots, \hat{J}_{1}$, by the formula \eqref{pCN_RPV_multi};

        \smallskip

        \item Compute $\phi^{\hat{j}} = \sum\limits_{\ell=1}^{L} \Phi^{\ell} \big( \tilde{\mathbf{Z}}_{n}^{\hat{j}}; \mathbf{Y}^{\ell} \big) $, $\hat{j}=1, 2, \cdots, \hat{J}_{1}$, with the model-data function $\Phi^{\ell}$ defined by the formula \eqref{model_data_misfit_l};

        \smallskip

        \item Sort the candidate sate $ \big\{ \tilde{\mathbf{Z}}_{n}^{\hat{j}} \big \}_{\hat{j}=1 }^{\hat{J}_{1}}$ corresponding the sequence $\big\{ \phi^{\hat{j}} \big\}_{\hat{j} =1}^{\hat{J}_{1}}$ from smallest to largest;

        \smallskip

        \item Set the candidate sate $ \tilde{\mathbf{Z}}_{n} $ by the first sate of the sorted candidate sates $ \big\{ \tilde{\mathbf{Z}}_{n}^{\hat{j}} \big \}_{\hat{j}=1 }^{\hat{J}_{1}}$.

    \end{enumerate}
\end{algorithm}

It is easy to find that the Algorithm \ref{Gibbs_Random_Proposal_multi} has improved Algorithm \ref{Gibbs_Random_Proposal} by using smaller integers $J_{0}, J_{1}$ and $J_{2}$. In particular, the number $J_{2}$ could be taken as $J_{2}=1$, since the state $\{ \mathbf{Z}_{j} \}_{j = J_{1}}^{J_{0}}$ are independent of each other. However, the computation cost of the Algorithm \ref{Gibbs_Random_Proposal_multi}, given by
\begin{equation}
    T^{\text{multi}} = T_{0}\hat{J}_{1}NJ_{0},
    \label{cost_Gibbs_multi}
\end{equation}
turns out to be greatly expensive.
The ratio between the computational cost of the Algorithms \ref{Gibbs_gPC} and \ref{Gibbs_Random_Proposal_multi} is
\begin{equation}
    R_{T} = \frac{ \hat{T} }{ T^{\text{multi}} }
    = \frac{ (T_{1} + T_{0}\hat{J}_{2})NJ_{0} }{ T_{0}\hat{J}_{1}NJ_{0} }
    = \frac{ T_{1} + T_{0}\hat{J}_{2}  }{ T_{0}\hat{J}_{1} }.
    \label{cost_ratio}
\end{equation}
Compared with the Algorithm \ref{Gibbs_Random_Proposal_multi},
 the Algorithm \ref{Gibbs_gPC} also explores $\hat{J}_{1}$ states in each iteration. However, using the surrogate model $\big\{ \tilde{\mathcal{G}}^{\ell} \big\}_{\ell=1}^{L}$,
  the Algorithm \ref{Gibbs_gPC} has an advantage that it gives a rough estimate of these $\hat{J}_{1}$ states with cheaper computational cost.
 Noting that $\hat{J}_{2}$ could be smaller if the surrogate model is more accurate, and the number $T_{1}$ could decrease if the total number of the basis function of the set $\big\{ \mathcal{H}_{\alpha}(\mathbf{Z}) \big\}_{\alpha \in \mathcal{I}, |\alpha| = 0}^{ |\alpha| = \tilde{N}}$ is smaller.
  In other words, we can reduce the computational cost with a smaller ratio $R_{T}$, if the surrogate model is more accurate and can be evaluated at cheaper computational cost. Further more, the computational cost $T_{0}$ will linearly increase with respect to  $L$ (the number of the incident waves) and $M$ (the number of the observation directions), while the computational cost $T_{1}$ is insensitive to these parameters. This implies that, in the case of a large number of incident waves and observation directions,
   the total computational cost can be significantly reduced by using the surrogate model.

\section{Numerical Examples}
\begin{figure}[htbp]
  \centering
  \includegraphics[width = 6in]{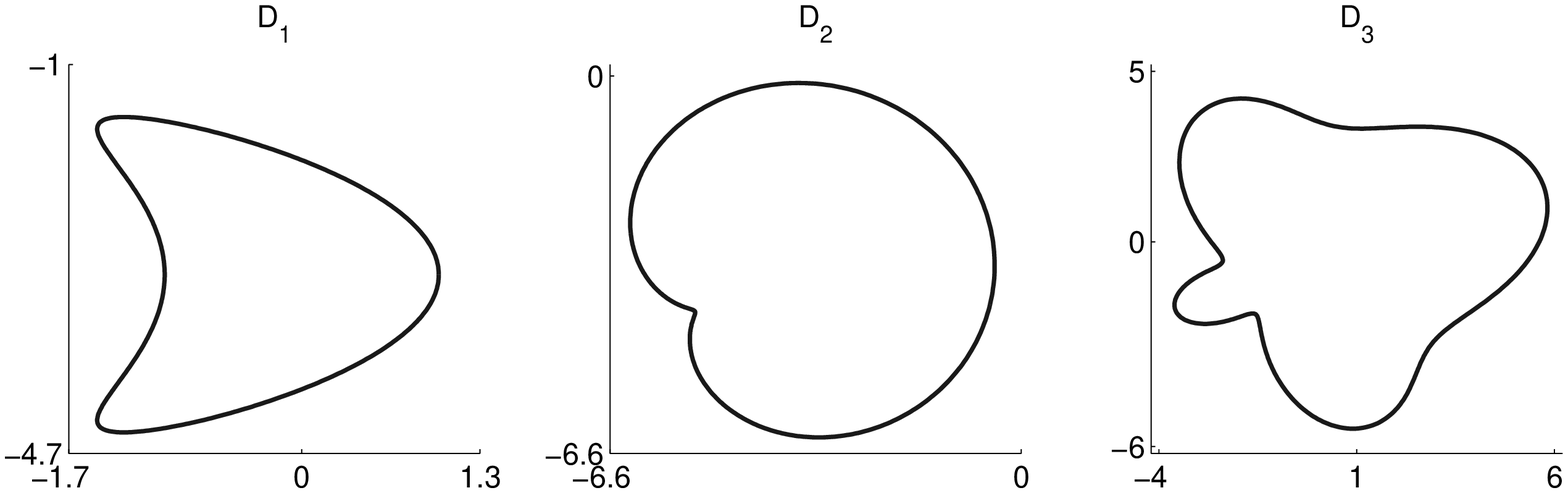}
  \caption{The boundaries of three acoustically sound-soft scatterers.}
  \label{three_scatter}
\end{figure}

In this section we exhibit numerical examples to demonstrate the effectiveness of the Bayesian method. In our former work \cite{Bayesian_phaseless_disk_line_2020_IP}, we consider inverse scattering of plane waves from a sound-soft disk with three unknown parameters, a line crack with four unknown parameters and a kite-shaped obstacle with six unknown parameters. In this paper, we shall take point source waves as incoming waves and extend the scenario to the kite-shaped obstacle with six unknowns in Example 1 and to other sound-soft obstacles with five (resp. eleven) unknowns in Example 2 (resp. Example 3).  In the final Example 4, we apply the gPC method to reduce the computation cost for recovering the kite-shaped obstacle. Below we represent the boundaries of three acoustically sound-soft scatterers in two dimensions (see Figure \ref{three_scatter}):
\begin{itemize}
    \item The first obstacle $D_{1}$,
        \begin{equation}
            \partial D_{1}
            = \left\{ x(t) = \big( - 0.65 + \cos t + 0.65 \cos (2t),
                                   \ -3 + 1.5 \sin t \big)^{\top},
                \quad 0\leq t \leq 2\pi \right\},
            \label{kite_t}
        \end{equation}
        which is the classical kite-shaped domain in inverse scattering problems.

        \smallskip

    \item The second obstacle $D_{2}$,
        \begin{eqnarray}
            \partial D_{2}
            &=& \left\{
                x(t) = \big( a + r(t)\cos t,
                             \ b + r(t) \sin t \big)^{\top},
                       \quad 0\leq t \leq 2\pi \right\},
                       \label{curve_t_type02_x_case01} \\
            r(t) &=& a_{0} + \sum_{n_{r}=1}^{N_{r}} a_{n_{r}}\cos(n_{r}t) + b_{n_{r}}\sin(n_{r}t),
            \quad 0\leq t \leq 2\pi, \label{curve_t_type02_r_case01}
        \end{eqnarray}
        with $N_{r} = 1$, $( a, b, a_{0}, a_{1}, b_{1} )^{\top} = (-5, -4, 2.5, 2, 1)^{\top}$. This obstacle can be regards as a disk with a notch.

        \smallskip

    \item The third obstacle $D_{3}$,
        \begin{eqnarray}
            \partial D_{3}
            &=& \left\{ x(t) = \big( a + r(t)\cos t,
                        \ b + r(t) \sin t \big)^{\top},
                \quad 0\leq t \leq 2\pi \right\},
                \label{curve_t_type02_x_case02} \\
            r(t) &=& a_{0} + \sum_{n_{r}=1}^{N_{r}} a_{n_{r}}\cos(n_{r}t) + b_{n_{r}}\sin(n_{r}t),
            \quad 0\leq t \leq 2\pi. \label{curve_t_type02_r_case02}
        \end{eqnarray}
        Here $N_{r} = 4$, $( a, b, a_{0}, a_{1}, b_{1}, \cdots, a_{4}, b_{4} )^{\top} = (-1, -1, 4, 2, 1, 0, 0, 0, 0, 0, 1)^{\top}$. There are 11 parameters in this scatterer. It will be
         used to investigate the effectiveness of the numerical method for recovering complex scatterers with high dimensional parameters.
\end{itemize}
\smallskip

We choose the Hausdorff distance (HD) to quantify the numerical error between the reconstructed and exact boundaries. The Hausdorff distance between two obstacles $\partial D_{1}$ and $\partial D_{2}$ is defined by
\begin{equation}
    d_{H}(\partial D_{1}, \partial D_{2})
    := \max\bigg\{ \sup_{x_{1} \in \partial D_{1}}
                   \inf_{x_{2}\in \partial D_{2}}
                     |x_{1} - x_{2}|, \quad
                   \sup_{x_{2} \in \partial D_{2}}
                   \inf_{x_{1} \in \partial D_{1}}
                   |x_{2} - x_{1}| \bigg\}.
    \label{Hausdorff_distance}
\end{equation}
The percent relative error (PRE) can be used to measure to what extent we have made use of the information of the observation data. In this paper, we denote by
 $f_{_{PRE}}$ the percent relative error (PRE) between the phaseless far-field data of the reconstructed boundary and the polluted observation. Notice that the polluted observation (\ref{observation_nosie_l}) is of the form of a $M\times L$ matrix, that is, $\mathbb{Y} = (\mathbf{Y}^{1}, \mathbf{Y}^{2}, \cdots, \mathbf{Y}^{L} ) \in \mathbb{R}^{M}\times \mathbb{R}^{L} $. One can calculate the phaseless far-field data $\tilde{\mathbb{Y}} = (\tilde{\mathbf{Y}}^{1}, \tilde{\mathbf{Y}}^{2}, \cdots, \tilde{\mathbf{Y}}^{L} )$ of the numerical reconstruction $\tilde{\mathbf{Z}}$ by
\begin{equation}
    \tilde{\mathbf{Y}}^{\ell} = \mathcal{G}^{\ell}(\tilde{\mathbf{Z}}),
    \quad \ell = 1, 2, \cdots, L,
    \label{forward_reconstruction}
\end{equation}
where $\mathcal{G}^{\ell}$ denotes the forward operator mapping the space of obstacle parameters to the space of phaseless far-field pattern in the discrete sense (see \eqref{forward_observation_l}). Then the percent relative error$f_{_{PRE}}$ is given by
\begin{equation}
    f_{_{PRE}} := \| \mathbb{Y} - \tilde{\mathbb{Y}} \|/\| \mathbb{Y} \|.
    \label{PRE_ff_phaseless}
\end{equation}
If $f_{_{PRE}}$ is small enough, we can claim that we have make full use of the observation data.

Recalling the definition \eqref{observation_nosie_l} of the polluted observation, we can construct two types of observations. Firstly, we consider an ideal setting where the observation data is the exact forward solution polluted by a special sample of noise pollution.
The accuracy of our numerical method can be examined in
 such an idea setting, that is,
the Hausdorff distance between the reconstructed and exact boundaries  indicates whether the numerical method is accurate or not.
Secondly, we consider a more practical setting where the exact forward solution is artificially polluted by a general sample of noise. In this case we can test the robustness of the numerical method in practical applications.

Unless otherwise stated, we always perform our numerical examples with the following assumptions.

\begin{itemize}
    \item The wave number is taken as $k=2$;

    \smallskip

    \item The incident waves are excited at source positions located at
        \begin{equation}
        x^{\text{in}}_{\ell} = (R \cos\theta_{\ell}, \ R \sin\theta_{\ell})^{\top},
        \quad \theta_{\ell} = 2\pi (\ell-1)/L,
        \quad \ell=1, 2, \cdots,L;
        \label{position_spherical}
    \end{equation}

    \item The observation directions are
    \begin{equation}
        \hat{\mathbf{x}}_{m} = \big(\cos\theta_{m}, \ \sin\theta_{m} \big)^{\top},
        \quad \theta_{m} = -\pi + 2\pi (m-1)/M,
        \quad m=1,2,\cdots,M;
        \label{direction_observation}
    \end{equation}

    \item For simplification, we only consider $L = M$;

    \smallskip

    \item In the setting of the prior distribution $P_{pr}$, the initial guess is supposed to be a unit circle centered at the origin with the variance variable $\sigma_{pr} =1$;

    \smallskip

    \item Corresponding to the $l$-th incident wave $u^{\text{in}}_{\ell}(x)$, the observation pollution $\eta^{\ell}$ $=$ $( \eta^{\ell}_{1}$, $\eta^{\ell}_{2}$, $\cdots$, $\eta^{\ell}_{M} )^{\top}$ is supposed to be a M-dimensional Gaussian variable, given by
        \begin{equation}
                \eta^{\ell}_{m}
                = \sigma_{\eta} \times
                  \left| u^{\infty}(\hat{\mathbf{x}}_m; \hat{\mathbf{Z}}, x^{\text{in}}_{\ell}, k) \right|
                  \ \omega^{\ell}_{m},
                  \quad m = 1, 2, \cdots, M,
                  \label{observation_noise_02}
          \end{equation}
        where $\hat{\mathbf{Z}}$ represents the exact obstacle parameters, $\omega^{\ell}_{m} \sim \mathcal{N}(0, 1)$, $m = 1, 2, \cdots, M$, and $\sigma_{\eta}$ is the noise coefficient. In other words, for $\ell = 1, 2, \cdots, L$, we take $\eta^{\ell} \sim \mathcal{N}( \mathbf{0}, \mathbf{\Sigma_{\eta^{\ell}}} )$ and the diagonal matrix $\mathbf{\Sigma_{\eta^{\ell}}} = \mbox{diag}( \sigma^{\ell}_{1},  \sigma^{\ell}_{2}, \cdots, \sigma^{\ell}_{M} )$ with $\sigma^{\ell}_{m} = \big( \sigma_{\eta}\times \left| u^{\infty}(\hat{\mathbf{x}}_m; \mathbf{Z}, x^{\text{in}}_{\ell}, k)\right| \big)^{2}, m = 1, 2, \cdots, M$. In our numerical tests, we choose $\sigma_{\eta} = 3\%, 6\%, 9\% $;

    \smallskip

    \item The observation $\mathbf{Y}^{\ell} = (y_{1}^{\ell}, y_{2}^{\ell}, \cdots, y_{M}^{\ell})^{\top}$, $\ell = 1, 2, \cdots, L$, is constructed as
        \begin{equation}
            y_{m}^{\ell}
            = \left| u^{\infty}(\hat{\mathbf{x}}_m; \hat{\mathbf{Z}}, x^{\text{in}}_{\ell}, k) \right|
              + \eta^{\ell}_{m}, \quad m = 1, 2, \cdots, M.
            \label{observation_noise_numerical}
        \end{equation}
        Here $\hat{\mathbf{Z}}$ is the exact obstacle parameters, and $\eta^{\ell}_{m}$, $\ell = 1, 2, \cdots, L$, $m = 1, 2, \cdots, M$ is the observation noise specified in \eqref{observation_noise_02};

    \smallskip

    \item In the ideal setting, the special sample of the observation noise is supposed to be $\omega^{\ell}_{m} \sim \mathcal{N}(0, 1)$ with $\omega^{\ell}_{m} = 0$, $\ell = 1, 2, \cdots, L$, $m = 1, 2, \cdots, M$, and $ \sigma_{\eta} = 3\%$;

    \smallskip

    \item In the practical settings of the observations, we generate 1000 samples of the observation noise $\omega^{\ell}_{m} \sim \mathcal{N}(0, 1)$, $\ell = 1, 2, \cdots, L$, $m = 1, 2, \cdots, M$. This leads to 1000 samples of the noise-polluted observation data with the noise coefficient $ \sigma_{\eta}$ and 1000 numerical reconstructions for these samples.
        We perform statistical analysis of these reconstructions to demonstrate the robustness of our numerical method;

    \smallskip

    \item The integral equation method is used for getting forward solutions. For this purpose we adopt the MATLAB code given by \cite[Chapter 8]{GNakamura_RPotthast_2015_InverseModeling};

    \smallskip

    \item All calculations are performed using MATLAB R2014a on a personal laptop with a 2.29 GHz CPU and 7.90 GB RAM.

\end{itemize}

\subsection{Example 1: Kite-shaped domain}

In this subsection, the kite-shaped sound-soft obstacle $D_{1}$ expressed by \eqref{kite_t} is supposed to be illuminated by incident waves \eqref{incident_wave_spherical}. The source positions are given by \eqref{position_spherical} with $R = 6$.
Obviously, the boundary of $D_{1}$ can be parameterized by six parameters
\begin{equation}
    \mathbf{Z}:= ( z_{1}, z_{2}, \cdots, z_{6} )^{\top}.
    \label{parameter_kite}
\end{equation}
The exact parameters are $\hat{\mathbf{Z}}$ $= (\hat{z}_{1}$, $\hat{z}_{2}$, $\cdots$, $\hat{z}_{6})^{\top} = (-0.65, -3, 1, 0.65, 1.5, 0)^{\top}$.
To implement the Algorithm \ref{Gibbs_Random_Proposal}, we choose $J_{0} = 20000$, $J_{1}=10000$, $J_{2}=100$ and $J_{3}=101$. The mean of the prior distribution $P_{pr}$ is $\mathbf{m}_{pr} = (0, 1, 0, 0, 1, 0)^{\top}$, as the initial guess is assumed to be a unit circle centered at the origin.

\begin{figure}[htbp]
  \centering
  \includegraphics[width = 5in]{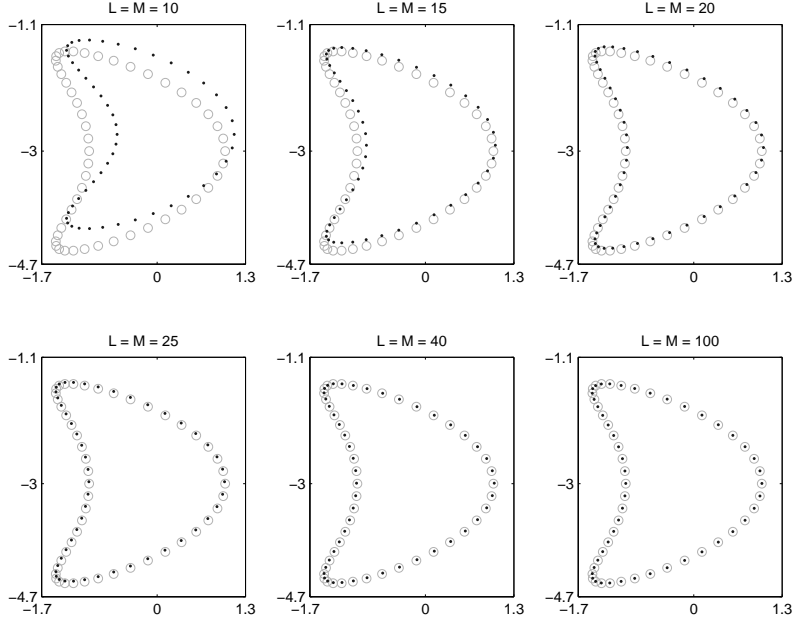}
  \caption{Reconstructions of the kite-shaped obstacle at the idea setting with different $L (M=L)$. The circle $\tch{\circ}$ and dot $\tcn{\cdot}$ denote the exact and reconstructed boundaries, respectively.}
  \label{kite_spherical_NoNoise_diff_NInOb}
\end{figure}

\begin{figure}[htbp]
  \centering
  \includegraphics[width = 3in]{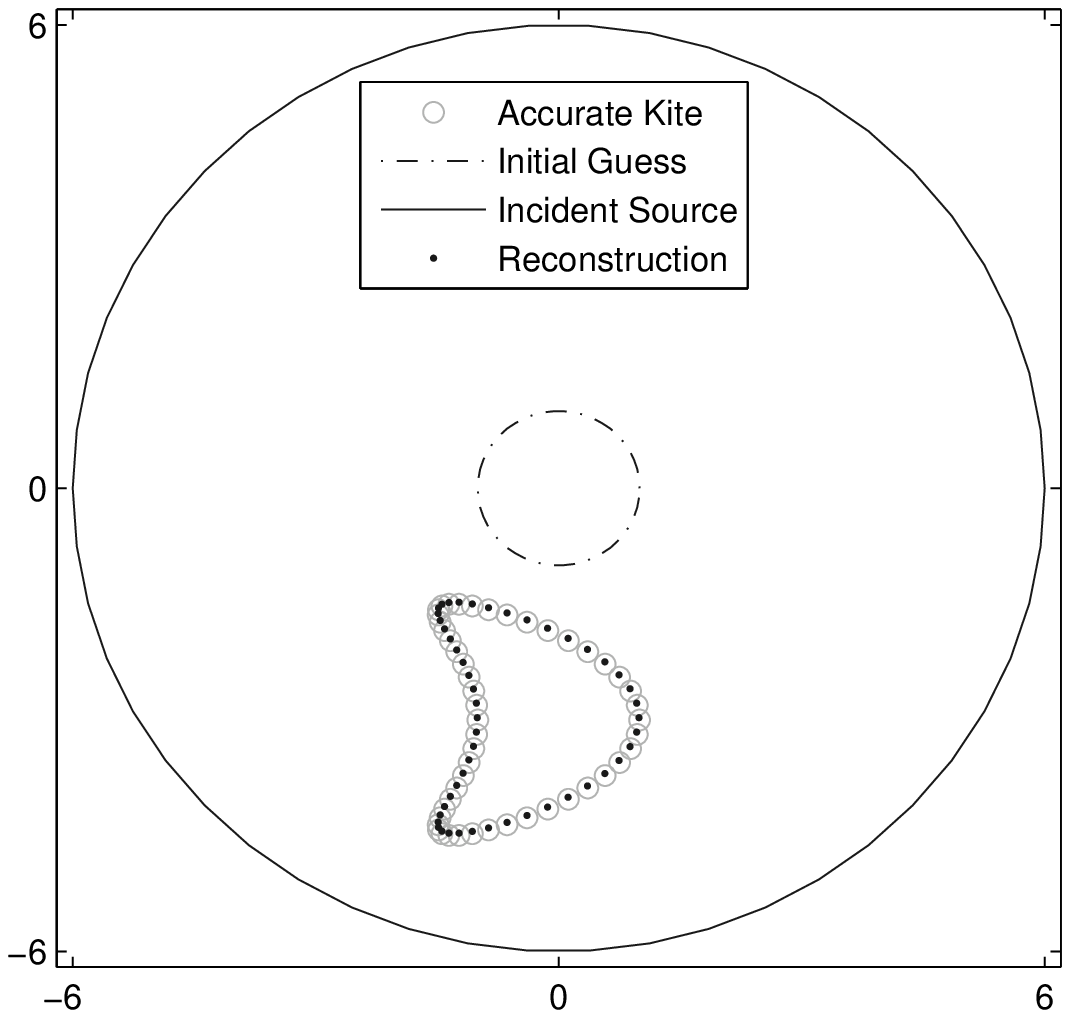}
  \caption{Reconstruction of a kite-shaped obstacle at the idea setting with $L=M=25$.}
  \label{kite_spherical_NoNoise_NInOb25}
\end{figure}

At first, we use ideal observations to investigate the accuracy of our method. For different choices of $L$ $(M=L)$, we exhibit in Table \ref{kite_spherical_NoNoise_diff_NInOb_table} and Figure \ref{kite_spherical_NoNoise_diff_NInOb} the numerical reconstructions and the Hausdorff distance between the numerical and  exact boundaries. We find that the reconstructed parameters are getting more accurate as the number of incident and observation directions becomes larger. The Hausdorff distances are less than $0.01$ if we choose $L$ and $M$ large enough such as $L = M = 40, 100$. The numerical solutions with $L = M = 15, 20$ are less accurate, since the Hausdorff distances are larger than 0.04 in these cases. The numerical solution with $L = M = 10$ turns out to be unreliable; see Figure \ref{kite_spherical_NoNoise_diff_NInOb}.
In Figure \ref{kite_spherical_NoNoise_NInOb25}, we plot
the initial guess,  positions of the incoming point sources together with the exact and reconstructed boundaries. Our numerical examples show that the inversion scheme
 does not depend on choice of initial guess.

\begin{table}[htbp]
    \centering
    \caption{Reconstruction of parameters and Hausdorff distance vs $L(M=L)$.}
    \label{kite_spherical_NoNoise_diff_NInOb_table}
    \begin{tabular}{c|c|c}
    \hline
      $L$ & Reconstructed parameters $z_j$ & HD  \\
    \hline
      10 & -0.4673,   -2.7488,    0.8614,    0.7405,    1.3993,    0.1154 & 0.1660  \\
    \hline
      15 & -0.6073,   -2.9121,    0.9478,    0.6881,    1.4661,    0.0431 & 0.0431  \\
    \hline
      20 & -0.6179,   -2.9475,    1.0013,    0.6423,    1.5169,   -0.0028 & 0.0489  \\
    \hline
      25 & -0.6526,   -2.9720,    0.9983,    0.6462,    1.4966,    0.0038 & 0.0216 \\
    \hline
      40 & -0.6479,   -2.9935,    1.0123,    0.6473,    1.4966,   -0.0028 & 0.0048 \\
    \hline
      100 & -0.6507,   -2.9975,    1.0010,    0.6520,    1.4988,    0.0018 & 0.0007 \\
    \hline
    \end{tabular}
\end{table}

\begin{figure}[htbp]
  \centering
  \includegraphics[width = 5.8in]{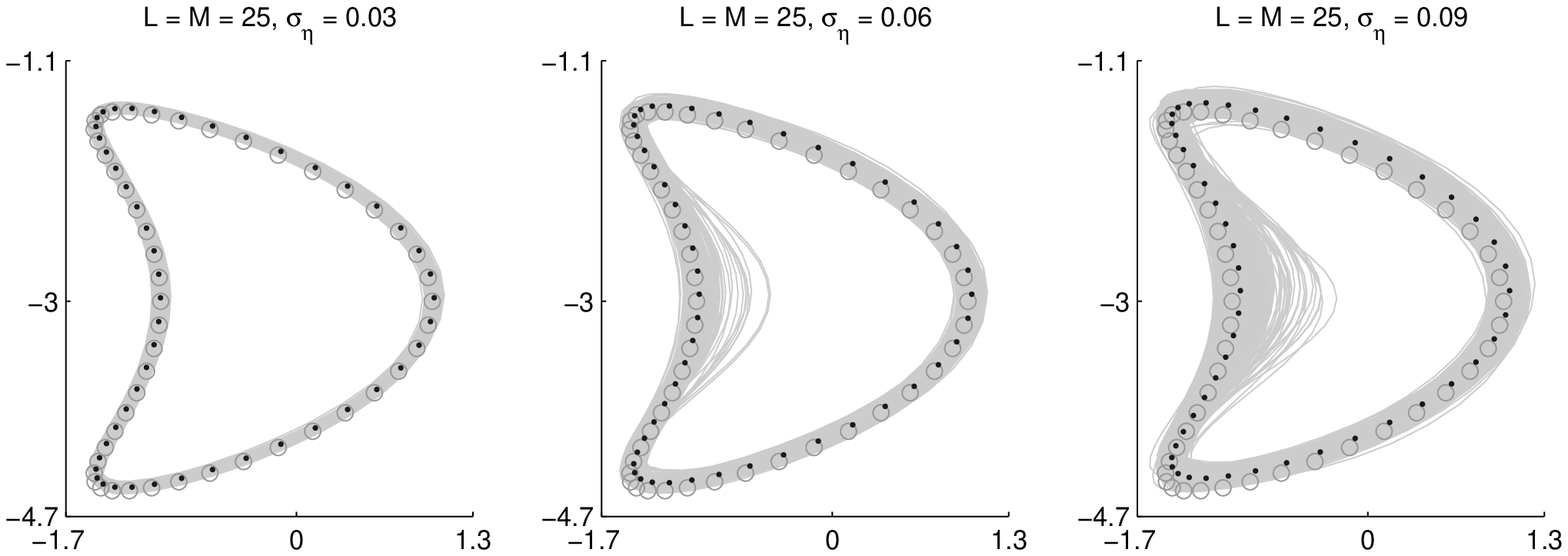}
  \caption{Reconstructions of the kite-shaped obstacle with $L=M=25$ at different noise levels $\sigma_{\eta} = 3\%$ (left), $6\%$ (center), $9\%$ (right). The circle $\tcn{\circ}$ denotes the exact boundary, the closed curves $\tcs{\overline{\quad}}$ are numerical reconstructions and the dot $\tcn{\cdot}$ is the mean of 1000 numerical reconstructions from noisy observations.}
  \label{kite_spherical_AddNoise_3_NInOb25}
\end{figure}

\begin{figure}[htbp]
  \centering
  \includegraphics[width = 5.8in]{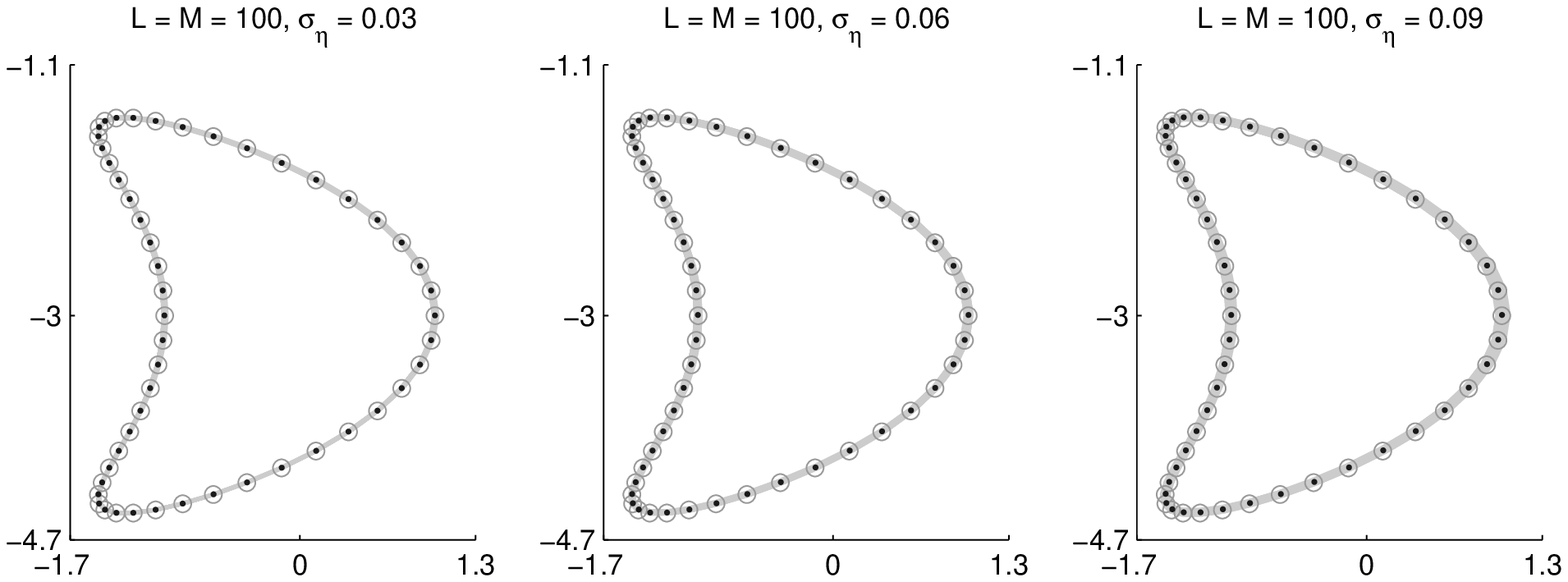}
  \caption{Reconstructions of the kite-shaped obstacle with $L=M=100$ at different noise levels $\sigma_{\eta} = 3\%$ (left), $6\%$ (center), $9\%$ (right). The circle $\tcn{\circ}$ is the exact boundary, the closed curves $\tcs{\overline{\quad}}$ are the numerical reconstructions and the dot $\tcn{\cdot}$ is the mean of 1000 numerical reconstructions from polluted observations.}
  \label{kite_spherical_AddNoise_3_NInOb100}
\end{figure}

Secondly, we consider the noise-polluted phaseless far-field data. The mean and standard deviation (SD) are used to perform statistical analysis over the 1000 reconstructions.
The mean of these reconstructions is reported in Table \ref{kite_spherical_AddNoise_3_NInOb25_table}, and
the mean and standard deviation (SD) of the Hausdorff distances at different noise levels $\sigma_{\eta}$ with $L=M=25$ are described in Figure \ref{kite_spherical_AddNoise_3_NInOb25}.
From Figure \ref{kite_spherical_AddNoise_3_NInOb25}, we find that the numerical reconstructions are distorted for large observation noise.
However, the results in Table \ref{kite_spherical_AddNoise_3_NInOb25_table} show that, our  approach is robust against the noise pollution, because the mean and standard deviation of the Hausdorff distances are relatively small even at the noise level of $9\%$.
Further, one can observe that less noise could give rise to  more reliable reconstructions.
 We also show the numerical results with $L=M=100$ in Figure \ref{kite_spherical_AddNoise_3_NInOb100}. From there one can claim that more reliable results can be obtained if we increase the number of observation directions in practical applications. Notice that, the numerical reconstruction is still accurate even if at the noise level of $\sigma_{\eta} = 0.09$. 

\begin{table}[htbp]
    \centering
    \caption{Numerical solutions vs $\sigma_{\eta}$ with $L=M=25$.}
    \label{kite_spherical_AddNoise_3_NInOb25_table}
    \begin{tabular}{c|c|c|c}
    \hline
    \  $\sigma_{\eta} $ & mean of reconstructed parameters $z_j$  & mean of HD & SD of HD \\
    \hline
     $3\%$ & -0.6372,   -2.9724,    1.0089,    0.6436,    1.4988,    0.0024  &  0.0229 &     0.0128 \\
    \hline
     $6\%$ &  -0.6215,   -2.9442,    1.0040,    0.6463,    1.4902,    0.0116  &  0.0467 &     0.0233 \\
    \hline
     $9\%$ &  -0.6038,   -2.9150,    0.9925,    0.6559,    1.4822,    0.0276   &   0.0659 &    0.0356 \\
    \hline
    \end{tabular}
\end{table}

\begin{figure}[htbp]
  \centering
  \includegraphics[width = 5.8in]{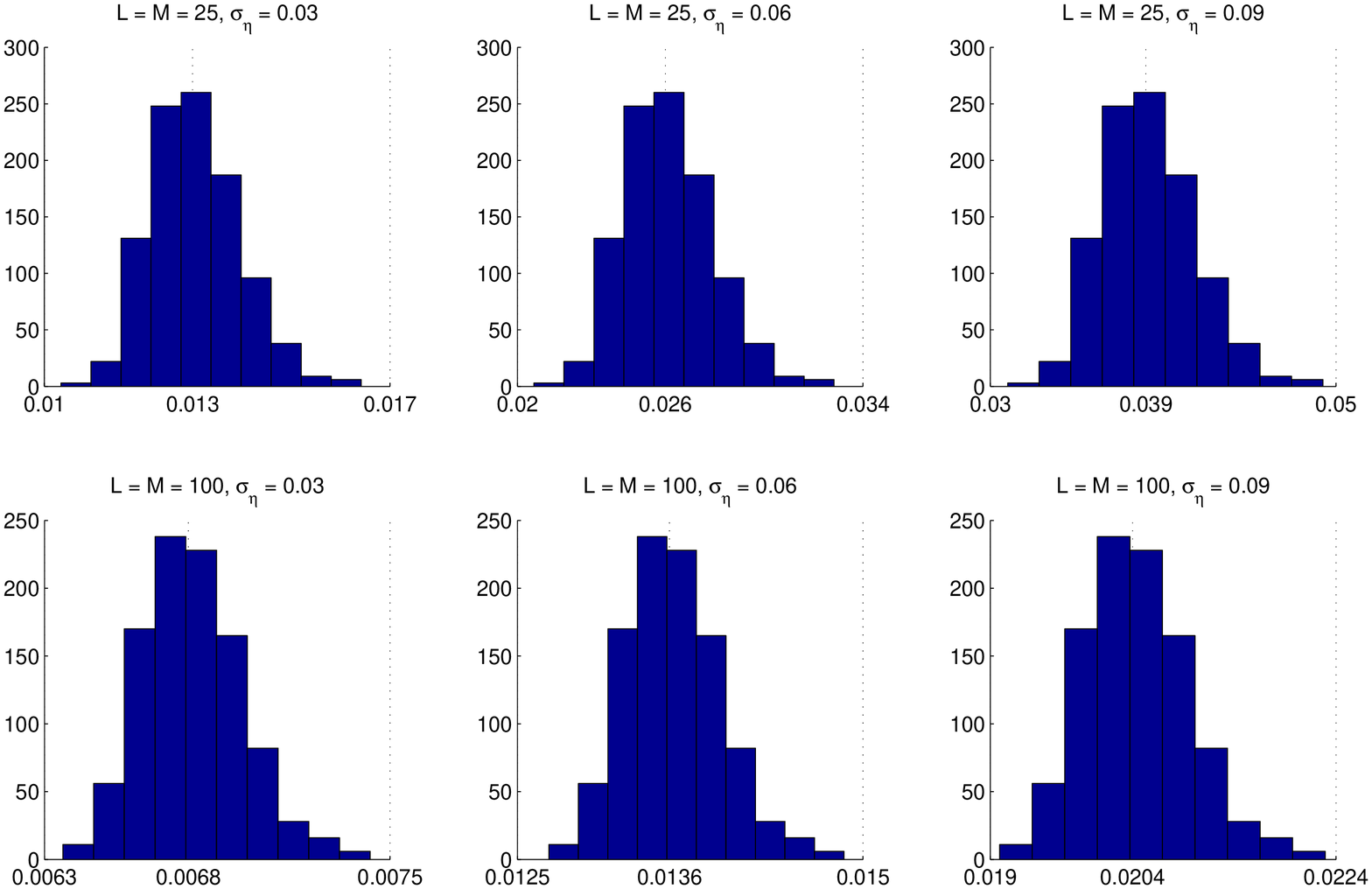}
  \caption{Histograms of $f_{_{PRE}}$ between the noise polluted observation and the exact phaseless data with $L=M=25$ (top), $L=M=100$ (bottom) at $\sigma_{\eta} = 3\%$ (left), $6\%$ (middle), $9\%$ (right).}
  \label{kite_spherical_ff_PRE_AddNoise_NInOb25_NInOb100}
\end{figure}

\begin{figure}[htbp]
  \centering
  \includegraphics[width = 5.8in]{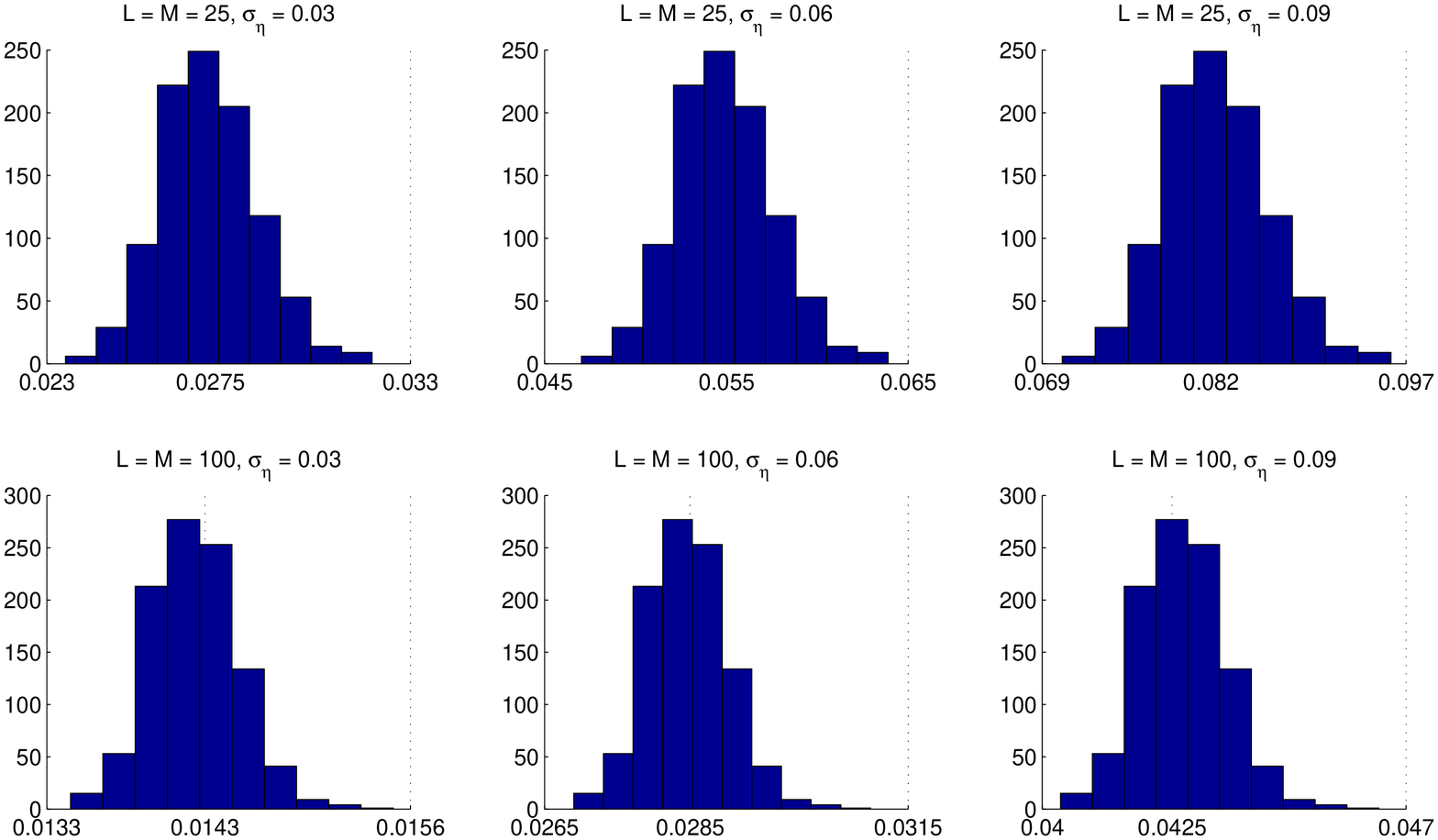}
  \caption{Histograms of $f_{_{PRE}}$ between the noise polluted observation defined by (\ref{observation_noise_type02}) and the exact phaseless far-field data with $L=M=25$ (top), $L=M=100$(bottom) at $\sigma_{\eta} = 3\%$ (left), $6\%$ (middle), $9\%$ (right).}
  \label{kite_spherical_ff_PRE_AddNoise_NInOb25_NInOb100_02}
\end{figure}

Using \eqref{PRE_ff_phaseless},  we can compute the percent relative error $f_{_{PRE}}$ between the noise polluted observation and the exact observation. Recall
from \eqref{observation_noise_02} that the phaseless far-field data are polluted by observation noise at the level $\sigma_{\eta}$. Let $\hat{\mathbb{Y}} = (\hat{\mathbf{Y}}^{1}$, $\hat{\mathbf{Y}}^{2}$, $\cdots$, $\hat{\mathbf{Y}}^{L} )$ be the exact phaseless far-field pattern of the obstacle $\hat{\mathbf{Z}}$, and let $\mathbb{Y}_{j} = (\mathbf{Y}^{1}_{j}$, $\mathbf{Y}^{2}_{j}$, $\cdots$, $\mathbf{Y}^{L}_{j} )$ be the 1000 samples of the noise polluted observation. Given $M=L$ and the noise level $\sigma_{\eta}$, we calculate the percent relative error $f^{j}_{_{PRE}} = \| \hat{\mathbb{Y}} - \mathbb{Y}_{j} \| / \| \hat{\mathbb{Y}} \|$ for $j=1$, 2, $\cdots$, 1000, and show their histograms in Figure \ref{kite_spherical_ff_PRE_AddNoise_NInOb25_NInOb100}. We observe that $\big \{ f^{j}_{_{PRE}} \big \}_{j=1}^{1000} $  decrease when the number of observation and incident directions increases, and they are extraordinary small in comparison with the noise level $\sigma_{\eta}$. Then we consider a second type of observation noise defined by
\begin{equation}
    \eta^{\ell}_{m}
    = \sigma_{\eta} \times
      \ \omega^{\ell}_{m},
      \quad \omega^{\ell}_{m} \sim \mathcal{N}(0, 1),
      \quad \ell = 1, 2, \cdots, L,
      \quad m = 1, 2, \cdots, M.
      \label{observation_noise_type02}
\end{equation}
The corresponding histograms of $\big \{ f^{j}_{_{PRE}} \big \}_{j=1}^{1000} $ are plotted in Figure \ref{kite_spherical_ff_PRE_AddNoise_NInOb25_NInOb100_02}, which are proven larger than those corresponding to
the observation noise \eqref{observation_noise_02}. In the case $L = M = 25$, the percent relative error is closer to the given noise level $\sigma_{\eta}$. 

\subsection{Example 2}

The boundary of the second sound-soft obstacle \eqref{curve_t_type02_x_case01}-\eqref{curve_t_type02_r_case01} can be parameterized by five parameters
\begin{equation}
    \mathbf{Z}:= ( z_{1}, z_{2}, \cdots, z_{5} )^{\top}= ( a, b, a_{0}, a_{1}, b_{1} )^{\top},
    \label{parameter_cure_type02}
\end{equation}
where the vector $(a, b)^{\top}$ indicates the position/location of this obstacle.
The exact obstacle parameters are $\hat{\mathbf{Z}}$ $= (\hat{z}_{1}$, $\hat{z}_{2}$, $\cdots$, $\hat{z}_{5})^{\top} = (-5, -4, 2.5, 2, 1)^{\top}$. In this example we shall test the capability of our approach for recovering a micro-subboundary (that is, the notch of $\partial D_2$) of a sound-soft obstacle.

We generate incident waves \eqref{incident_wave_spherical} by setting $R = 9$ in \eqref{position_spherical}. Let the mean of the prior distribution $P_{pr}$ be $\mathbf{m}_{pr} = (0, 0, 1, 0, 0)^{\top}$. In the Algorithm \ref{Gibbs_Random_Proposal}, we choose $J_{0} = 20000$, $J_{1}=10000$, $J_{2}=100$ and $J_{3}=101$.

\begin{figure}[htbp]
  \centering
  \includegraphics[width = 6in]{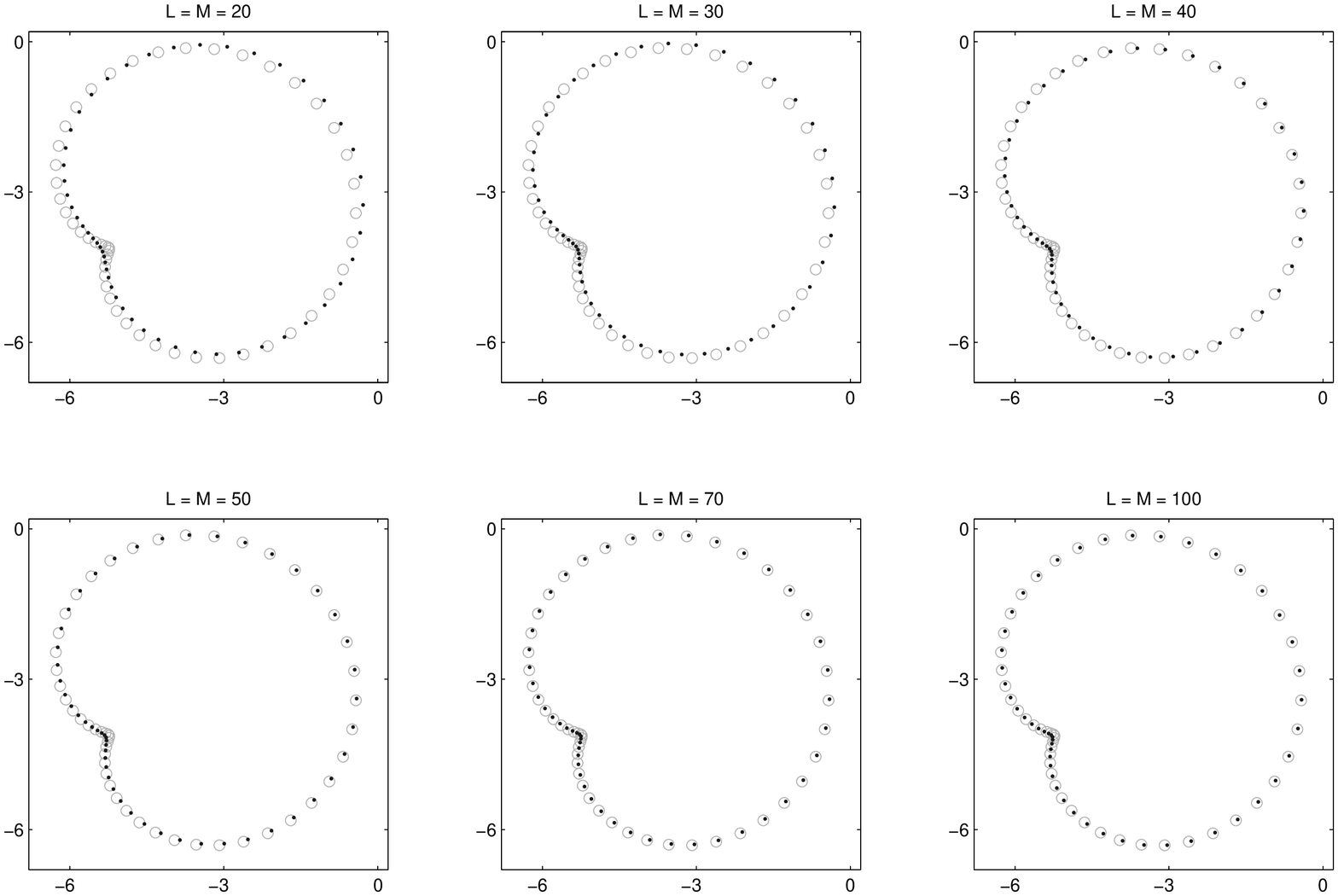}
  \caption{Recovery of the second obstacle at the idea setting with different $L$ $(M=L)$. The circle $\tcn{\circ}$ denotes the exact boundary and the dot $\tcn{\cdot}$ the reconstructions.}
  \label{Parameter02_Case01_spherical_NoNoise_diff_NInOb}
\end{figure}

\begin{figure}[htbp]
  \centering
  \includegraphics[width = 3in]{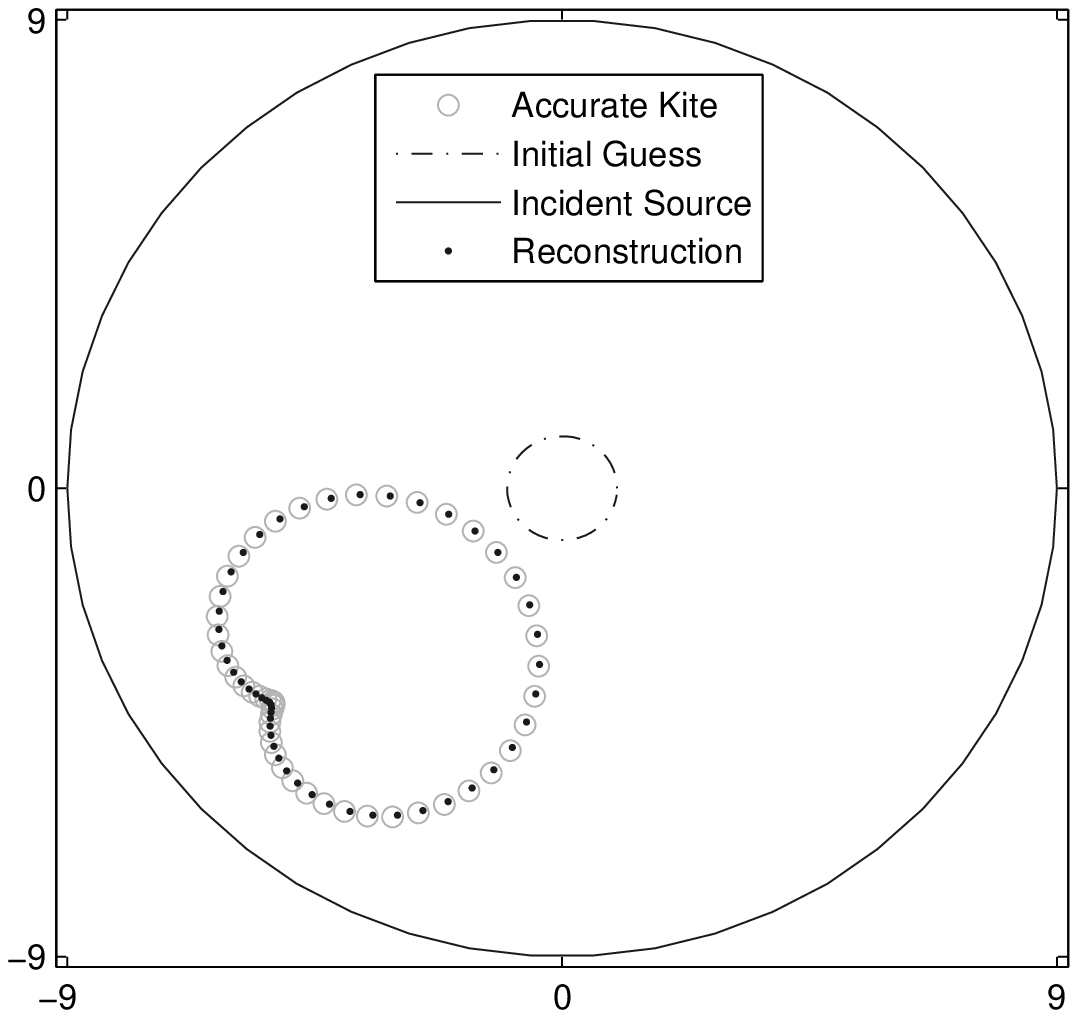}
  \caption{Reconstructions of the second obstacle at the idea setting with $L=M=50$.}
  \label{Parameter02_Case01_spherical_NoNoise_NInOb50}
\end{figure}

In the case of ideal observations, we discuss the accuracy of the numerical solutions for different choices of $L$ $(M=L)$ in Table \ref{Parameter02_Case01_spherical_NoNoise_diff_NInOb_table} and Figure \ref{Parameter02_Case01_spherical_NoNoise_diff_NInOb}. Obviously, the Hausdorff distance between the reconstructed and exact boundaries decreases as the number of incident and observation directions becomes larger. A satisfactory reconstruction of the location/position of the obstacle can be achieved even if $L, M$ are small such as $L = M = 20$. However, an accurate recovery of the notch requires large $L$ and $M$ such as $L = M = 50, 70, 100$.
Figure \ref{Parameter02_Case01_spherical_NoNoise_NInOb50}  shows that our approach is not sensitive to
 the initial guess.

\begin{table}[htbp]
    \centering
    \caption{Numerical Reconstructions vs $L(M=L)$.}
    \label{Parameter02_Case01_spherical_NoNoise_diff_NInOb_table}
    \begin{tabular}{c|c|c}
    \hline
      $L$ & Reconstructed parameters $z_j$ & HD \\
    \hline
      20 & -4.6808,   -3.8144,    2.6498,    1.6899,    0.8117 & 0.1434 \\
    \hline
      30 &  -4.7925,   -3.8674,    2.6087,    1.8064,    0.9058 & 0.0763 \\
    \hline
      40 &  -4.8550,   -3.9437,    2.5645,    1.8497,    0.9180 & 0.0500 \\
    \hline
      50 &  -4.9163,   -3.9550,    2.5437,    1.8923,    0.9479 & 0.0194 \\
    \hline
      70 &   -4.9581,   -3.9744,    2.5227,    1.9565,    0.9842 & 0.0187 \\
    \hline
      100 &   -4.9603,   -3.9846,    2.5254,    1.9466,    0.9696 & 0.0071 \\
    \hline
    \end{tabular}
\end{table}

\begin{figure}[htbp]
  \centering
  \includegraphics[width = 6in]{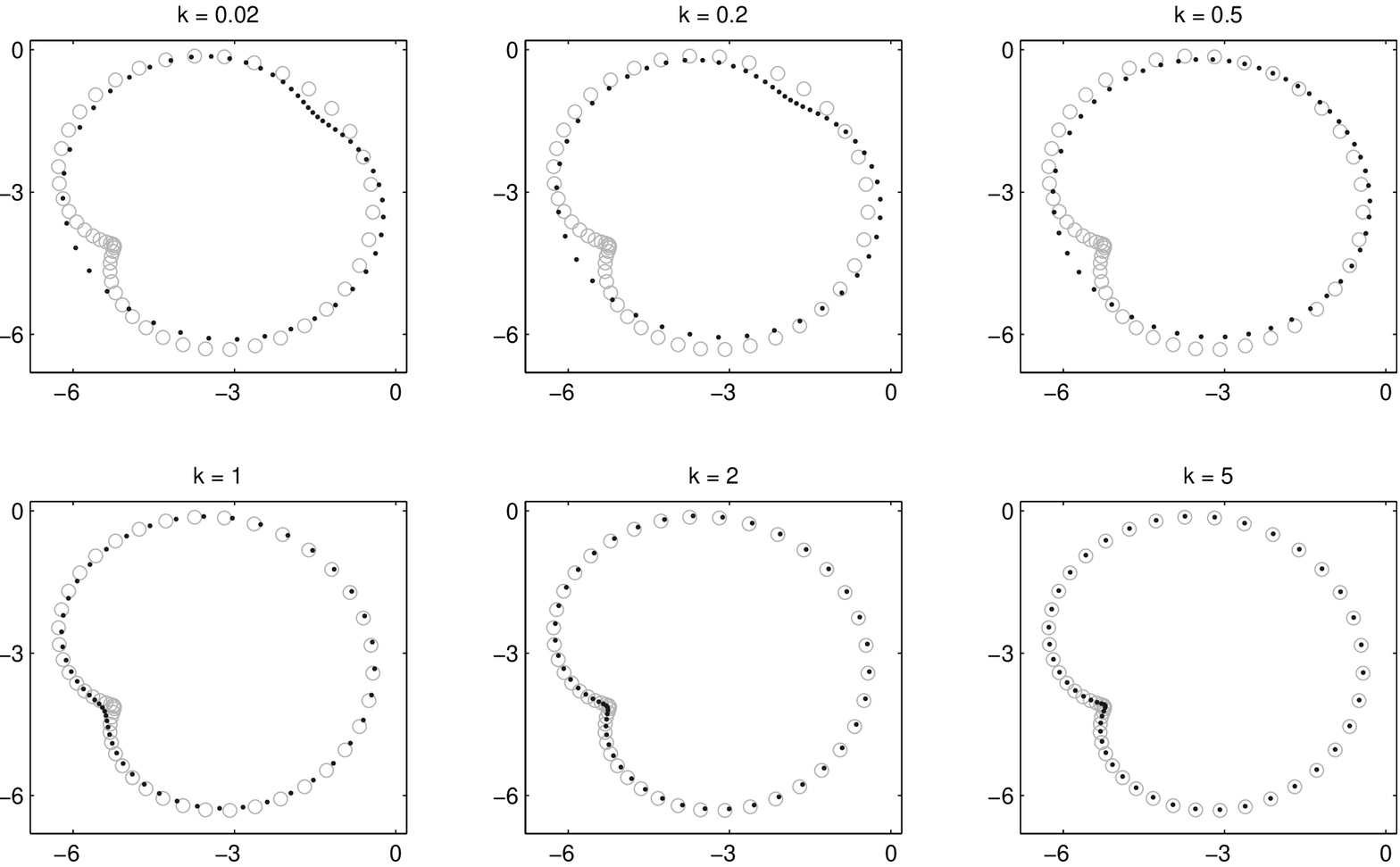}
  \caption{Reconstructions of the second obstacle with different wave numbers at the idea setting. We set $L = M = 50$. The circle $\tcn{\circ}$ and dot $\tcn{\cdot}$
  denote respectively the exact and reconstructed boundaries.}
  \label{parameter02_case01_spherical_NoNoise_diff_k}
\end{figure}

Now we discuss the accuracy of the numerical solutions at different wave numbers.
Since the wave number $k$ is inversely proportional to the wave length, it is more difficult to extract the information of the notch at smaller wavenumbers.
Hence, numerical reconstructions cannot be expected to be accurate when the wave number $k$ becomes smaller; see the results shown in Table \ref{parameter02_case01_spherical_NoNoise_diff_k_table} and Figure \ref{parameter02_case01_spherical_NoNoise_diff_k} with $L=M=50$ in the idea setting.
The Hausdorff distance decreases as the wave number $k$ increases. We also calculate the percent relative error $f_{_{PRE}}$ between the phaseless far-field data of the reconstructed and polluted observations in Table \ref{parameter02_case01_spherical_NoNoise_diff_k_table}. The decreasing $f_{_{PRE}}$ for larger wavenumbers is consistent with the reconstructed Hausdorff distance.
For the chosen wavenumbers, see Figure \ref{parameter02_case01_spherical_NoNoise_diff_k}, one can always get a good approximation of obstacle location, while an accurate reconstruction of the notch can be achieved only for large $k$.

\begin{table}[htbp]
    \centering
    \caption{The percent relative error $f_{_{PRE}}$ and numerical reconstructions vs wave number $k$ with $L=M=50$.}
    \label{parameter02_case01_spherical_NoNoise_diff_k_table}
    \begin{tabular}{c|c|c|c}
    \hline
      $k$ & $f_{_{PRE}}$ & Reconstructed parameters $z_j$ & HD \\
    \hline
      0.02 & $0.37\%$ & -2.1874,   -2.1006,    2.6857,   -1.1903,   -1.1858 & 0.3109  \\
    \hline
      0.2 & $0.67\%$ & -2.4166,   -1.9278,    2.6838,   -0.9258,   -1.3725 & 0.3126 \\
    \hline
      0.5 & $1.36\%$ & -2.7712,   -2.5510,    2.8827,   -0.4930,   -0.5941 & 0.1612 \\
    \hline
      1 & $0.43\%$ & -4.7948,   -3.8822,    2.6234,    1.7193,    0.8458  & 0.0399 \\
    \hline
      2 & $0.25\%$ &  -4.9163,   -3.9550,    2.5437,    1.8923,    0.9479  & 0.0194 \\
    \hline
      5 & $0.16\%$ & -4.9984,   -3.9879,    2.4955,    2.0064,    1.0079 & 0.0105 \\
    \hline
    \end{tabular}
\end{table}

\begin{figure}[htbp]
  \centering
  \includegraphics[width = 6in]{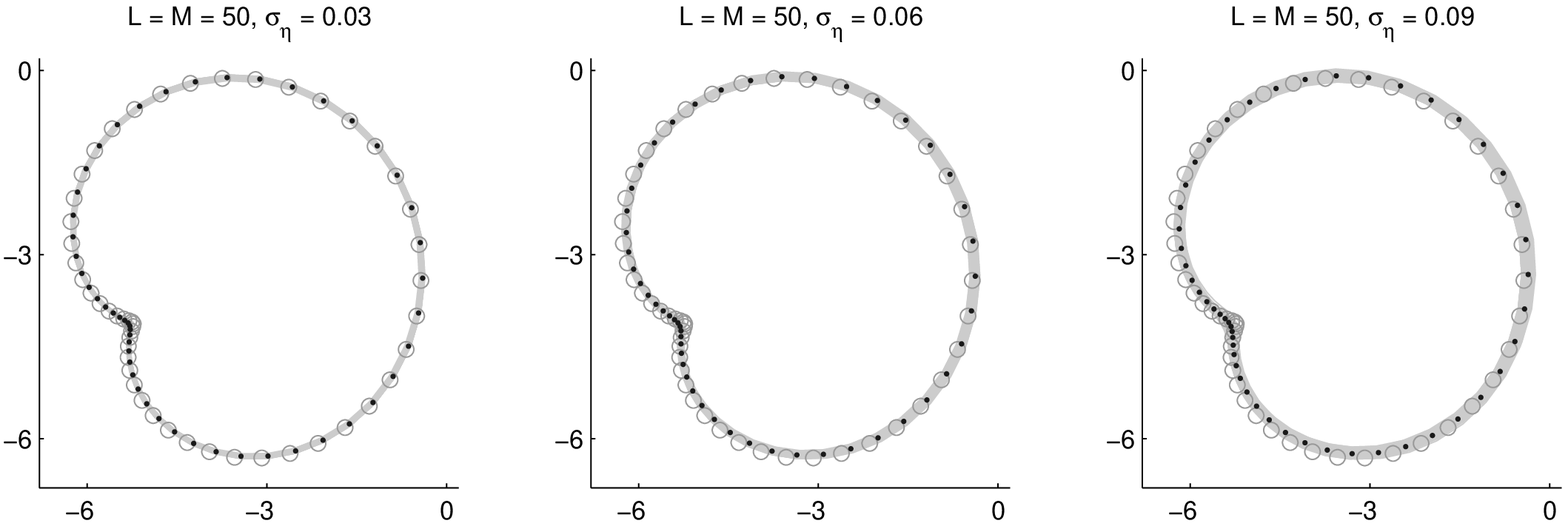}
  \caption{Reconstructions of the second obstacle with $L=M=50$ at different noise levels $\sigma_{\eta} = 3\%$ (left), $6\%$ (center), $9\%$ (right). The circle $\tcn{\circ}$ is the exact boundary, the closed curves $\tcs{\overline{\quad}}$ are numerical reconstructions and the dot $\tcn{\cdot}$ is the mean of 1000 numerical reconstructions from polluted observation data.}
  \label{Parameter02_Case01_spherical_AddNoise_3_NInOb50}
\end{figure}

Using polluted data, we get numerical reconstructions
 for 1000 samples of the observation noise. 
 In Table \ref{Parameter02_Case01_spherical_AddNoise_3_NInOb50_table}
 we calculate the mean of each reconstructed parameter, the mean and standard deviation (SD) of the Hausdorff distances at different noise levels $\sigma_{\eta}$ with $L=M=50$.
It can be observed from Figure \ref{Parameter02_Case01_spherical_AddNoise_3_NInOb50} that the reconstructed boundaries
maybe inaccurate if the noise level is large, because the reconstruction of the notch is blurred when the noise level increases. But the numerical method is robust against the noise pollution, as the mean and standard deviation of the Hausdorff distance are small; see Table \ref{Parameter02_Case01_spherical_AddNoise_3_NInOb50_table}.

\begin{table}[htbp]
    \centering
    \caption{Numerical solutions vs $\sigma_{\eta}$ with $L=M=50$.}
    \label{Parameter02_Case01_spherical_AddNoise_3_NInOb50_table}
    \begin{tabular}{c|c|c|c}
    \hline
    \  $\sigma_{\eta} $ & mean of reconstructed parameters  & mean of HD & SD of HD \\
    \hline
     $3\%$ &  -4.9116,   -3.9528,    2.5479,    1.8949,    0.9482   &  0.0301 &     0.0135 \\
    \hline
     $6\%$ &  -4.8488,   -3.9173,    2.5739,    1.8322,    0.9177  &  0.0541 &     0.0199 \\
    \hline
     $9\%$ &   -4.7936,   -3.8844,    2.5926,    1.7835,    0.8938  &   0.0702 &    0.0249 \\
    \hline
    \end{tabular}
\end{table}

\subsection{Example 3}

The boundary of the third obstacle \eqref{curve_t_type02_x_case02}-\eqref{curve_t_type02_r_case02} can be parameterized by the following $11$ parameters
\begin{equation}
    \mathbf{Z}:= ( z_{1}, z_{2}, \cdots, z_{} )^{\top}= ( a, b, a_{0}, a_{1}, b_{1}, \cdots, a_{4}, b_{4} )^{\top},
    \label{parameter_cure_03}
\end{equation}
and the exact obstacle parameters are $\hat{\mathbf{Z}}$ $= (\hat{z}_{1}$, $\hat{z}_{2}$, $\cdots$, $\hat{z}_{11})^{\top} = (-1, -1, 4, 2, 1, 0, 0, 0, 0, 0, 1)^{\top}$.

We generate incident point source wave \eqref{incident_wave_spherical} with $R = 8$ in \eqref{position_spherical}. The mean of the prior distribution $P_{pr}$ is set to be $\mathbf{m}_{pr} = (0, 0, 1, 0, \cdots, 0)^{\top}$. To implement the Algorithm \ref{Gibbs_Random_Proposal}, we choose $J_{0} = 50000$, $J_{1}=40000$, $J_{2}=100$ and $J_{3}=101$.

\begin{figure}[htbp]
  \centering
  \includegraphics[width = 5.8in]{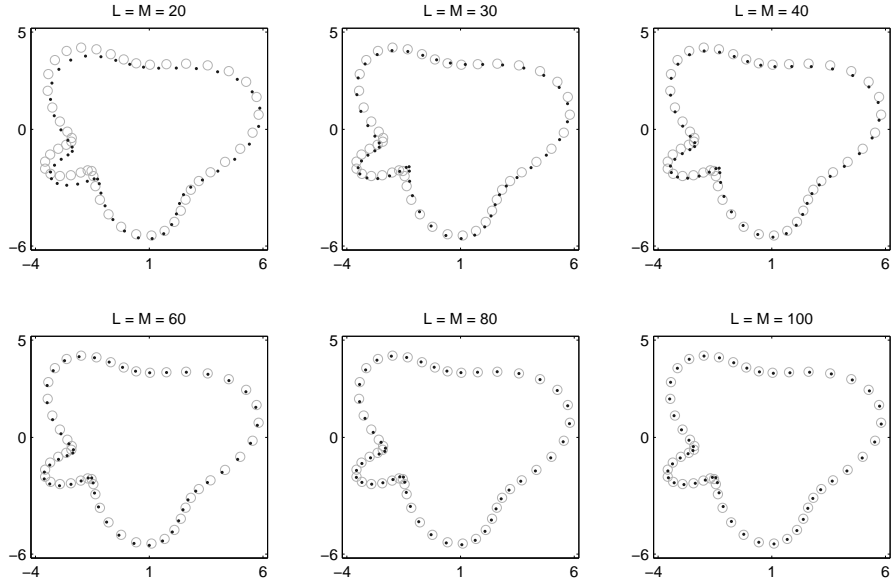}
  \caption{Reconstructions of the third obstacle at the idea setting with different $L$ $(M=L)$. The circle $\tcn{\circ}$ denotes the exact boundary and the dot $\tcn{\cdot}$ denotes the reconstructed boundary.}
  \label{Parameter02_Case04_spherical_NoNoise_diff_NInOb}
\end{figure}

\begin{figure}[htbp]
  \centering
  \includegraphics[width = 6in]{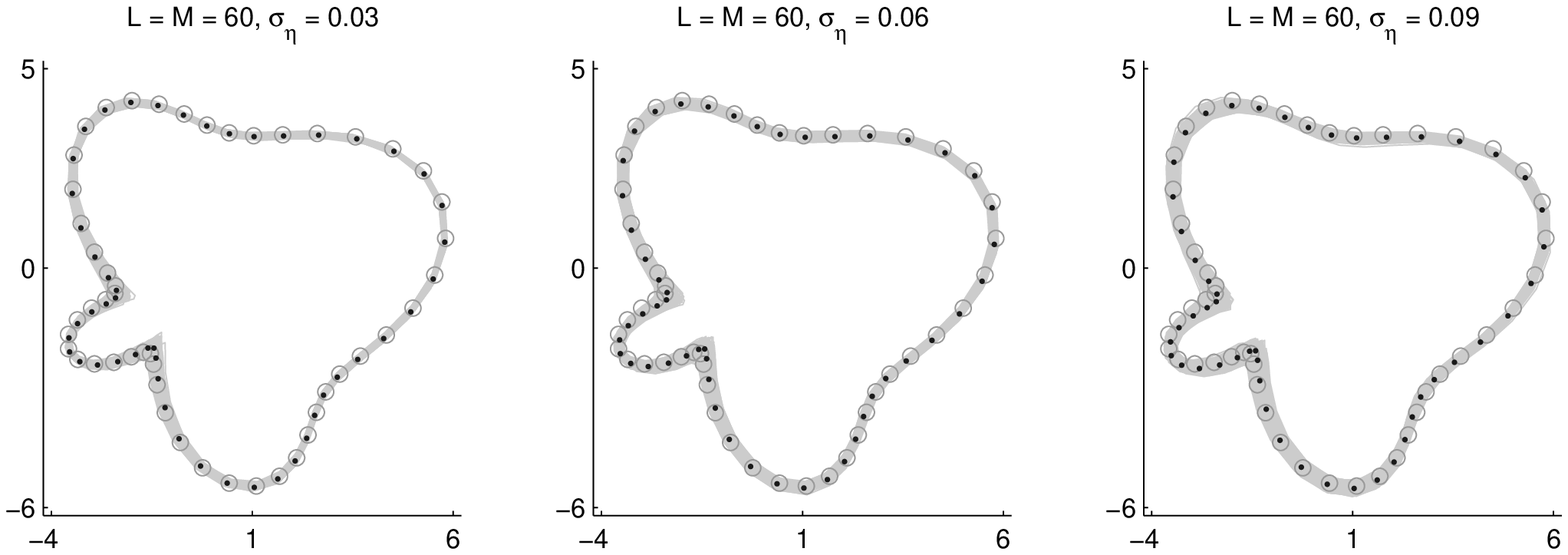}
  \caption{Reconstructions of the third obstacle with $L=M=60$ at different noise levels $\sigma_{\eta} = 3\%$ (left), $6\%$ (center), $9\%$ (right). The circle $\tcn{\circ}$ is the exact boundary, the closed curves $\tcs{\overline{\quad}}$ are numerical reconstructions and the dot $\tcn{\cdot}$ is the mean of 100 numerical reconstructions from polluted data.}
  \label{Parameter02_Case04_spherical_AddNoise_3_NInOb60}
\end{figure}

In the noise-free case, the accuracy of the numerical solutions for different choice of $L (M=L)$ is shown in Table \ref{Parameter02_Case04_spherical_NoNoise_diff_NInOb_table} and Figure \ref{Parameter02_Case04_spherical_NoNoise_diff_NInOb}. Although there are 11 unknown parameters, an accurate reconstruction can be obtained if we have enough input and output information such as $L = M = 60, 80, 100$.
In the noisy case, as we have done in the previous examples, the mean and standard deviation are again used to statistically analyze 1000 reconstructions, which corresponds to 1000 samples of the observation noise; see Table \ref{Parameter02_Case04_spherical_AddNoise_3_NInOb60_table_01} and Figure \ref{Parameter02_Case04_spherical_AddNoise_3_NInOb60}. As observed in the previous examples, the numerical reconstruction maybe inaccurate if the observation noise is large in Figure \ref{Parameter02_Case04_spherical_AddNoise_3_NInOb60}. However, the numerical method is still robust against the noise pollution, because the mean and standard deviation of  the Hausdorff distance are small.

\begin{table}[htbp]
    \centering
    \caption{Numerical reconstructions  vs $L(M=L)$.}
    \label{Parameter02_Case04_spherical_NoNoise_diff_NInOb_table}
    \begin{tabular}{c|c|c}
    \hline
      $L$ & Reconstructed parameters $z_j$ & HD \\
    \hline
      20 & -0.790,   -1.543,    3.923,    1.821,    1.317,   -0.079,    0.148,   -0.089,    0.094,   -0.033,  0.961  & 0.1495  \\
    \hline
      30 & -0.986,   -1.329,    3.941,    2.003,    1.365,    0.027,   -0.063,   -0.166,   -0.009,   -0.047,  0.997  & 0.0688 \\
    \hline
      40 & -0.966,   -1.251,    3.904,    2.014,    1.224,    0.020,   -0.039,   -0.090,   -0.035,   -0.028,  1.003  & 0.0453 \\
    \hline
      60 & -1.075,   -1.131,    3.962,    2.070,    1.099,   -0.021,    0.037,   -0.046,   -0.044,   -0.002,  1.031  & 0.0100 \\
    \hline
      80 & -1.047,   -1.051,    3.963,    2.084,    1.042,   -0.007,    0.014,   -0.036,   -0.039,    0.005,  1.028  & 0.0074 \\
    \hline
      100 & -0.980,  -1.058,    3.988,    1.981,    1.080,    0.022,   -0.023,   -0.027,    0.005,  -0.010,   1.011  & 0.0039 \\
    \hline
    \end{tabular}
\end{table}

\begin{table}[htbp]
    \centering
    \caption{Numerical solutions vs $\sigma_{\eta}$ with $L=M=60$.}
    \label{Parameter02_Case04_spherical_AddNoise_3_NInOb60_table_01}
    \begin{tabular}{c|c|c|c}
    \hline
    \  $\sigma_{\eta} $ & mean of reconstructed parameters $z_j$ & mean of HD & SD of HD\\
    \hline
     $3\%$ & -1.02,   -1.10,   3.96,    2.04,    1.11,    0.02,   -0.01,   -0.05,   -0.01,    1.03   &  0.0233 & 0.0161  \\
    \hline
     $6\%$ & -1.02,   -1.14,   3.94,    2.04,    1.14,    0.01,    0.002,  -0.06,   -0.01,    1.03   &  0.0319 & 0.0213  \\
    \hline
     $9\%$ & -1.00,   -1.20,   3.94,    2.01,    1.17,   -0.0004,  0.01,   -0.06,   -0.02,    1.02   &  0.0401 & 0.0256  \\
    \hline
    \end{tabular}
\end{table}

\subsection{Example 4: gPC method for the kite-shaped domain}

To demonstrate the efficiency of the surrogate model constructed by the gPC method, we apply the Algorithm \ref{Gibbs_gPC} to reconstruct the kite-shaped domain $D_{1}$, and compare the results with those obtained in Example 1.

Recall the exact obstacle parameters $\hat{\mathbf{Z}}$ $= (\hat{z}_{1}$, $\hat{z}_{2}$, $\cdots$, $\hat{z}_{6})^{\top} = (-0.65$, $1, 0.65, -3, 1.5, 0)^{\top}$. As in previous settings we set $R = 6$ and let the mean of the prior distribution $P_{pr}$ be $\mathbf{m}_{pr} = (0, 1, 0, 0, 1, 0)^{\top}$. For the surrogate model, we take $\tilde{N}=9$. To perform the Algorithm \ref{Gibbs_gPC}, we choose $J_{0} = 100$, $J_{1}=50$, $J_{2}=1$, $J_{3}=51$, $\hat{J}_{1}=1000$ and $\hat{J}_{2}=100$.

First, we need to calculate the chaos coefficients $u_{\alpha}^{\ell, m}$, $\alpha \in \mathcal{I}, |\alpha| = 0, 1, \cdots, \tilde{N}$, $\ell = 1, 2, \dots, L$, $m = 1, 2, \dots, M$. For this purpose, 6000 samples of the obstacle parameters $\mathbf{Z}$ are generated from the prior distribution $P_{pr}( \mathbf{Z} )$. Then, for every given $L$ and $M$, the integration in the equation \eqref{u_alpha_l_m} is calculated through the Monte Carlo method \cite{Robert_2004_MC}. The CPU time spend on calculating the chaos coefficients with $L = M = 25$ is 1.57 minutes.

\begin{figure}[htbp]
  \centering
  \includegraphics[width = 5.4in]{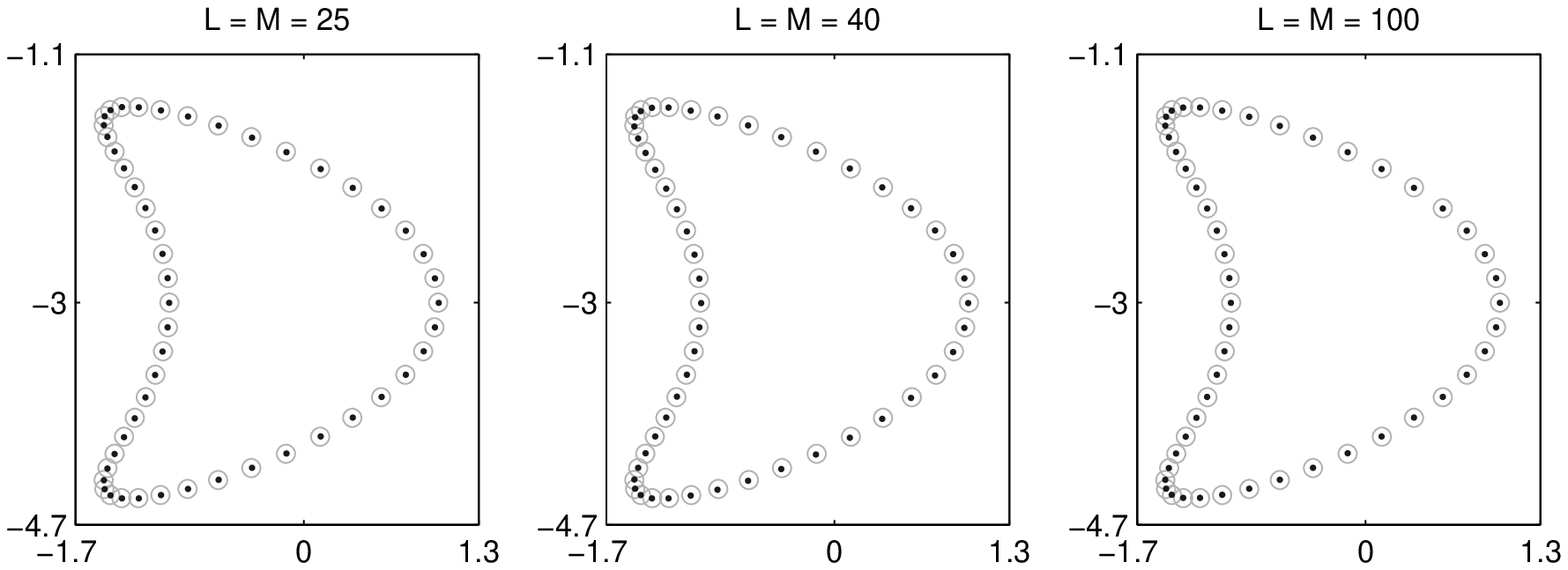}
  \caption{Reconstructions of the kite-shaped obstacle by the Algorithm \ref{Gibbs_gPC} at the idea setting with different $L$ $(M=L)$. The circle $\tcn{\circ}$ and dot $\tcn{\cdot}$ denote respectively
  the exact and reconstructed boundaries.}
  \label{gPC_kite_spherical_NoNoise_diff_NInOb}
\end{figure}

The numerical results with $L = M = 25, 40, 100$ are shown in Table \ref{gPC_kite_spherical_NoNoise_diff_NInOb_table} and Figure \ref{gPC_kite_spherical_NoNoise_diff_NInOb} in the noise-free case.
We observe that the reconstructions using the surrogate model are more accurate than those reported in Example 1; cf. Table \ref{kite_spherical_NoNoise_diff_NInOb_table}.
This is due to the reason that the Algorithm \ref{Gibbs_gPC} in the current Example 4 has explored $100000$ states in the iteration of the MCMC method,  while the Algorithm \ref{Gibbs_Random_Proposal} in Example 1 only explores $20000$ states.

Then we implement the Algorithm \ref{Gibbs_Random_Proposal} at $L = M = 25$ by exploring $100000$ states, with the results shown in Table \ref{kite_spherical_NoNoise_25_table_compare_gPC}.
 Comparing the results in Tables  \ref{kite_spherical_NoNoise_25_table_compare_gPC} and \ref{kite_spherical_NoNoise_diff_NInOb_table}, we find that  more accurate reconstructions can be achieved by exploring more states, since the Hausdorff distance decreases from 0.0216 (see Table \ref{kite_spherical_NoNoise_diff_NInOb_table} with $L=M=25$ ) to 0.0163 (see Table \ref{kite_spherical_NoNoise_25_table_compare_gPC}). But the resulting Hausdorff distance 0.0163 is still much larger than
 the reconstructed distance 0.0019 using the surrogate model
 . On the other hand, the computational cost with the surrogate model (18.68 minutes) is much cheaper than that in the Algorithm \ref{Gibbs_Random_Proposal} (81.37 minutes). It is worthy noting that the CPU time of 18.68 minutes also includes the computational cost (1.57 minutes) for calculating the chaos coefficients. In summary, using the surrogate model constructed by the gPC method, we have indeed reduced the computational cost.

\begin{table}[htbp]
    \centering
    \caption{Numerical reconstructions vs $L$ $(M=L)$ by the Algorithm \ref{Gibbs_gPC}.}
    \label{gPC_kite_spherical_NoNoise_diff_NInOb_table}
    \begin{tabular}{c|c|c}
    \hline
      $L$ & Reconstruction of parameters $z_j$ & HD \\
    \hline
      25 &  -0.6479,   -3.0005,    1.0030,    0.6485,    1.4999,   -0.0025 & 0.0019  \\
    \hline
      40 &  -0.6514,   -3.0032,    0.9935,    0.6515,    1.4981,    0.0053 & 0.0016  \\
    \hline
      100 & -0.6506,   -3.0005,    1.0006,    0.6495,    1.4970,   -0.0005 & 0.0006  \\
    \hline
    \end{tabular}
\end{table}

\begin{table}[htbp]
    \centering
    \caption{Comparison of numerical reconstructions using
    Algorithms \ref{Gibbs_Random_Proposal} and \ref{Gibbs_gPC}
     with $L = M = 25$.}
    \label{kite_spherical_NoNoise_25_table_compare_gPC}
    \begin{tabular}{c|c|c|c}
    \hline
      method & Reconstructed parameters $z_j$ & HD & tic-toc time \\
    \hline
      Algorithm \ref{Gibbs_Random_Proposal} & -0.6467,   -2.9721,    1.0106,    0.6502,    1.4947,    0.0080  & 0.0163  & 81.37 minutes  \\
    \hline
      Algorithm \ref{Gibbs_gPC} &  -0.6479,   -3.0005,    1.0030,    0.6485,    1.4999,   -0.0025 & 0.0019  &  18.68 minutes \\
    \hline
    \end{tabular}
\end{table}

\begin{figure}[htbp]
  \centering
  \includegraphics[width = 6in]{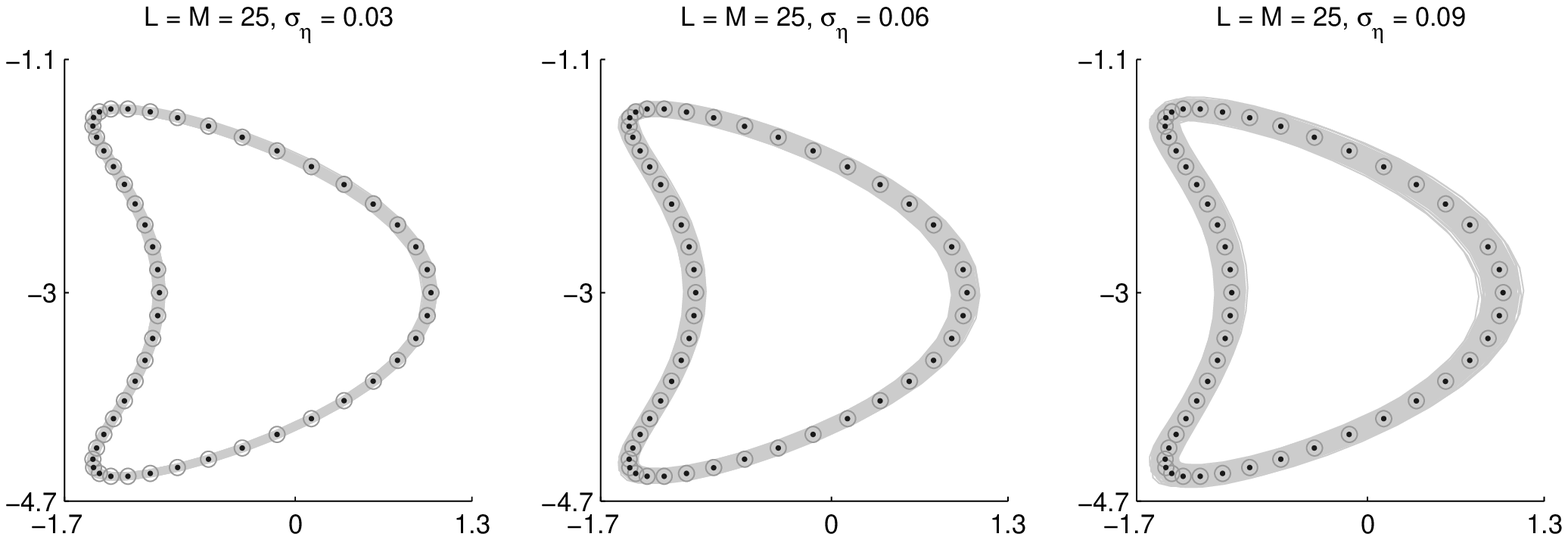}
  \caption{Reconstructions of the kite-shaped obstacle by the Algorithm \ref{Gibbs_gPC} with $L=M=25$ at different noise levels $\sigma_{\eta} = 3\%$ (left), $6\%$ (center), $9\%$ (right). The circle $\tcn{\circ}$ is the exact boundary, the closed curves $\tcs{\overline{\quad}}$ are the numerical reconstructions from 1000 samples of the observation noise and the dot $\tcn{\cdot}$ is the mean of these 1000 reconstructions.}
  \label{gPC_kite_spherical_AddNoise_3_NInOb25}
\end{figure}

Finally, the numerical results from the noisy data of the level $\sigma_{\eta}$ are shown in Table \ref{gPC_kite_spherical_AddNoise_3_NInOb25_table} and Figure \ref{gPC_kite_spherical_AddNoise_3_NInOb25}. As in Example 1, the numerical method is robust against the noise pollution, and the phaseless data with less noise give rise to a more reliable reconstruction result. Further more, with much cheaper computational cost, the  mean and standard deviation (SD) of Hausdorff distances in Table \ref{gPC_kite_spherical_AddNoise_3_NInOb25_table} are smaller than the corresponding results in Table \ref{kite_spherical_AddNoise_3_NInOb25_table}. Hence,  the surrogate model has  improved the numerical results reported in the first example.

\begin{table}[htbp]
    \centering
    \caption{Numerical solutions vs $\sigma_{\eta}$ with $L=M=25$ by the Algorithm \ref{Gibbs_gPC}.}
    \label{gPC_kite_spherical_AddNoise_3_NInOb25_table}
    \begin{tabular}{c|c|c|c}
    \hline
    \  $\sigma_{\eta} $ & mean of reconstructed parameters $z_j$  & mean of HD & SD of HD \\
    \hline
     $3\%$ & -0.6499,   -3.0000,    0.9999,    0.6497,    1.4999,  -0.0001  &  0.0076 & 0.0062 \\
    \hline
     $6\%$ & -0.6500,   -3.0003,    0.9991,    0.6497,    1.4997,  -0.0001  &  0.0158 & 0.0130 \\
    \hline
     $9\%$ & -0.6499,   -3.0003,    0.9988,    0.6496,    1.4997,  -0.0001  &  0.0238 & 0.0197 \\
    \hline
    \end{tabular}
\end{table}

\section{Conclusion}

In this paper, we apply the Bayesian approach to inverse time-harmonic scattering problems of recovering sound-soft obstacles from the phaseless far-field data excited by point source waves. 
Special attention has been paid to complex obstacles with high-dimensional parameters. When the dimension of unknown parameters becomes larger, both the computational cost of the forward model and the number of iteration steps in the MCMC method would increase sharply, giving rise to prohibitively high cost in the MCMC simulation.
We adopt the Gibbs method and the stochastic surrogate model based on the generalized polynomial chaos method to overcome this challenge. We also develop a strategy to combine the stochastic surrogate model with the MCMC method. In our numerical examples, the efficiency of these schemes are demonstrated without sacrificing too much accuracy.
Our future efforts will be devoted to recovering physical properties (e.g., refractive index) of an acoustically scatterer, which contain more parameters than its geometrical shape as discussed here. In this paper, the total number of boundary parameters is assumed to be known in advance. Removing or relaxing is assumption would also lead to
 a high dimensional inverse scattering problem,  and then more efficient schemes are needed in designing the MCMC method and the reduced model. Research outcomes along these directions will be reported in our forthcoming publications.

\bibliographystyle{plain}
\bibliography{ref_Bayesian_Acoustic_phaseless_curve}

\end{document}